%% file: 2004-21.tex
\documentclass{gtart_h}  

\input gtoutput

\lognumber{211}
\volumenumber{8}\papernumber{21}\volumeyear{2004}
\pagenumbers{779}{830}
\received{16 October 2001}
\revised{29 March 2004}
\accepted{29 April 2004}
\published{20 May 2004}
\proposed{Simon Donaldson}
\seconded{Ronald Stern, Robion Kirby}

\usepackage{amssymb, amsmath, amscd}

\numberwithin{equation}{section}
\newtheorem{thm}{Theorem}[section]
\newtheorem{lem}{Lemma}[section]
\newtheorem{defn}{Definition}[section]
\newtheorem{prop}{Proposition}[section]
\newtheorem{rem}{Remark}[section]
\newtheorem{exmp}{Example}[section]
\newtheorem{sublem}{Sublemma}[section]
\newtheorem{cor}{Corollary}[section]

\newcommand{\definition}[1]{\textit{#1}}
\newcommand{\Supp}{\operatorname{Supp}}
\newcommand{\Spec}{\operatorname{Spec}}
\newcommand{\Hom}{\operatorname{Hom}}
\newcommand{\AHS}{\operatorname{AHS}}
\newcommand{\Aut}{\operatorname{Aut}}
\newcommand{\Ad}{\operatorname{Ad}}
\newcommand{\id}{\operatorname{id}}
\newcommand{\const}{\mathit{const}}
\newcommand{\coker}{\operatorname{coker}}
\newcommand{\vol}{\mathit{vol}}
\newcommand{\im}{\operatorname{im}}
\newcommand{\interior}{\operatorname{int}}

\begin{document}

\title{ASD moduli spaces over four--manifolds\\with tree-like ends}
\asciititle{ASD moduli spaces over four-manifolds with tree-like ends}
\author{Tsuyoshi Kato}
\address{Department of Mathematics, Kyoto University\\Kyoto 606-8502, Japan} 
\email{tkato@math.kyoto-u.ac.jp}

\begin{abstract}
In this paper we construct Riemannian metrics and weight functions
over Casson handles. We show that the corresponding
Atiyah--Hitchin--Singer complexes are Fredholm for some class of Casson
handles of bounded type.  Using these, the Yang--Mills moduli spaces
are constructed as finite dimensional smooth manifolds over Casson
handles in the class.
\end{abstract}
\asciiabstract{%
In this paper we construct Riemannian metrics and weight functions
over Casson handles. We show that the corresponding
Atiyah-Hitchin-Singer complexes are Fredholm for some class of Casson
handles of bounded type.  Using these, the Yang-Mills moduli spaces
are constructed as finite dimensional smooth manifolds over Casson
handles in the class.}


\maketitlepage

\section*{Introduction}
\addcontentsline{toc}{section}{Introduction}
\subsection{Review of previous works on exotic smooth structures on
open four--manifolds}
Four--manifold theory has been deeply developed from two approaches.
One is based on geometry and functional analysis. In particular gauge
theory, Yang--Mills theory or Seiberg--Witten theory, construct moduli
spaces as the sets of solutions of nonlinear PDEs on four--manifolds.
They have been discovered to contain extremely rich information on
smooth structure.  In particular the construction of Donaldson's
invariant uses Yang--Mills moduli spaces. The invariant distinguishes
many mutually-homeomorphic but non-diffeomorphic pairs of smooth
four--manifolds.  The other is based on topology, in particular
Casson--Freedman's theory.  In high dimensional differential topology,
Whitney's trick to remove self-intersections of immersed discs has
played one of the most important r\^oles, but does not work in four
dimensions.  Casson's idea was, instead of removing them, to increase
self-intersections by attaching immersed two handles constructively
so that all self-intersections were able to be pushed away to infinity.
Freedman verified that these infinite towers made from immersed
two--handles were all homeomorphic to the standard open two--handle, and
they are called \definition{Casson handles}.  This allowed the complete
classification of oriented simply-connected topological four--manifolds
by their intersection forms and Kirby--Siebenmann classes
\cite{freedman}.

Combination of Casson--Freedman's theory with gauge theory provided
very deep results on smooth structures on four--manifolds.  One of the
most important results on exotic $\mathbb{R}^4$s was due to Taubes, who
discovered uncountably many exotic smooth structures on $\mathbb{R}^4$
\cite{taubes}.  In the essential step, he constructed Yang--Mills gauge theory
on periodic open four--manifolds.  Let $X$ be an oriented topological
four--manifold.  $X \backslash pt$  can admit at least one smooth
structure \cite{freedman and quinn}.  The idea was to verify that $X \backslash pt$ was
not able to admit any end-periodic smooth structure.  In fact otherwise
it would produce generic Yang--Mills moduli spaces over $X \backslash pt$
as  smooth manifolds.  On the other hand, detailed analysis on the
moduli spaces verified that it was impossible for such spaces to exist.
Thus the smooth structure of the end was sufficiently complicated
to obtain uncountably many exotic $\mathbb{R}^4$s.

\subsection{Casson handles}
Casson handles (CH) can be constructed inside
smooth four--manifolds.  One defines a smooth structure on CH by
restriction. In another words, a tower consists of smoothly-immersed
two--handles, and so  the tower itself admits  a smooth structure.
Each building block by an immersed two handle is called a
\definition{kinky handle}.
Even though any Casson handle is homeomorphic to the standard open
two--handle,  many of them are not diffeomorphic to the standard one.
In fact even for the simplest Casson handle,  CH$_+$,  consisting of a
single and positive kink at each stage, the following is known:
\begin{lem}{\rm\cite{bizaca2,bizaca and gompf}}\qua
$CH_+$ is exotic in the sense that the attaching circle 
cannot bound a smoothly embedded disc inside $CH_+$.
\end{lem}

For each signed infinite tree there corresponds a CH as a smooth
four--manifold.  In this paper we will define a subclass of signed
infinite trees which are called \definition{homogeneous trees of bounded
type}. These are constructed iteratively by attaching infinitely many
half-periodic trees. A connected infinite subtree of a homogeneous
tree of bounded type is called a \definition{tree of bounded type}.
Notice that a Casson handle of a tree of bounded type admits a smooth
embedding by another CH  of a tree of homogeneously bounded type.  In this
paper we will use the term \definition{Casson handle of bounded type} to
refer to any Casson handle constructed from a signed tree of homogeneously
bounded type.

Now any Casson handle of bounded type can be smoothly embedded in $CH_+$
preserving the attaching circles.  The above immediately implies the
following:
\begin{cor}
Any CH of bounded type is exotic.
\end{cor}

Typically, smooth structures on Casson handles will have a deep effect
on the smooth types of four--manifolds which contain them.  Thus
Casson handles are smooth open four--manifolds with the attaching
region, and comprise a rich class among open four--manifolds.  They
will provide highly nontrivial examples to study in open manifold
theory.  It is difficult to understand the important numerical
relationship between the growth of signed trees and the complexity of
smooth structures on the corresponding Casson handles.  In
\cite{gompf} some relation between Stein structures and the number of
kinks was found.

\subsection{Outline of the article}
In this paper we study Yang--Mills gauge theory over Casson handles.
Our final aim will be to measure the complexity of smooth structures
on Casson handles by means of gauge theory. Roughly speaking in order
to construct  Yang--Mills theory on open four--manifolds, one needs
to overcome two steps.  The first is the Fredholm theory of the linearized
equation. The other is perturbation theory, or transversality, where
the setting of perturbation  is  different from the closed case (as we
will explain below).  A general procedure in Yang--Mills theory tells
us that in a situation on four--manifolds where the Fredholm theory
is applicable, one can obtain Yang--Mills moduli spaces as generically
smooth manifolds of finite dimension.  In order to induce information
on smooth four--manifolds from these moduli spaces, they are required
to be non-empty under generic perturbation, where transversality theory
will enter effectively.

In this paper we show that for Casson handles of bounded type, one can
impose complete Riemannian metrics and weight functions on them so that
Fredholm theory can be done.  So we will obtain Yang--Mills moduli
spaces on them as generically smooth manifolds of finite dimension.
Perturbation theory on Casson handles is at present under development.

Let $(S,g)$ be a Riemannian (possibly open) oriented four--manifold, and
$E \to S$ be an $SO(3)$ bundle. A connection $A$ on $E$ is called
\definition{anti-self-dual} (\definition{ASD}) if its curvature form
satisfies the equation $F_A^+ \equiv  F_A+ *F_A =0$.  Roughly speaking
when $S$ is closed, the set of ASD connections modulo the gauge group is
the Yang--Mills moduli space, which is generically a smooth manifold of
finite dimension. In this case one does not have to take care so much
on the underlying function spaces.  Any choice among various Sobolev
spaces $W^k(S)$ ($k$ large) gives the same moduli space.

When $S$ is open, there are no standard choices of function spaces.
Let $w\co S \to [0, \infty)$ be a smooth function. Then one obtains
weighted Sobolev spaces $W^k_w(S)$. Let us choose an ASD connection
$A_0$, and consider two sets. One is $\widetilde{\mathfrak{M}}$,  the set of
ASD connections such that their curvature forms are in $L^2$. The other
is $\widetilde{\mathfrak{M}}(A_0) = \{A : A~\text{is ASD},~A - A_0 \in W^k_w\}$.
In general these spaces are very different and it seems difficult to
study the geometry of the former space.  In a periodic case,  one can see
that these give  the same moduli space modulo gauge transformation, which
follows from the exponential decay estimate on the curvature forms.
A standard argument on transversality also works using the decay
estimate.  For  the case of Casson handles, it seems  natural to regard
these  spaces $\widetilde{\mathfrak{M}}$ and $\widetilde{\mathfrak{M}}(A_0)$ as having
mutually different natures.   The exoticness of smooth types of various
Casson handles come from how the boundary solid tori are attached.  These
data are put in the cylindrical direction in our choice of metric.  So it
would be a delicate matter how these moduli spaces behave near the ends.
The complexity of these moduli spaces will reflect that of the Casson
handles.

In our class of Casson handles of bounded type,  we will show that
a Fredholm theory can be constructed using $\widetilde{\mathfrak{M}}(A_0)$.
Our main concern here is analysis of the Atiyah--Hitchin--Singer complex:
$$0 \longrightarrow W^{k+1}_w((Y,g)) \stackrel{d}{\longrightarrow}
W^k_w((Y,g);\Lambda^1) \stackrel{d^+}{\longrightarrow}
W^{k-1}_w((Y,g);\Lambda^2_+) \longrightarrow 0.$$
In this paper we will explicitly construct Riemannian metrics and weight
functions over Casson handles, so that the above complex is Fredholm
and its cohomology groups are calculable when restricted to our class.
This is the Fredholm part which we mentioned above.  One may generalize
this to the case with coefficients in the adjoint bundle $\Ad({\mathfrak{g}})$
of $E$.  We denote by $H^*_{A'}(\AHS)$ the cohomology groups with
coefficient at $A' \in {\mathfrak{M}}$.

Suppose the above complex is Fredholm. Then roughly speaking the ASD
moduli space $\widetilde{\mathfrak{M}}(A_0)/{\mathfrak{G}}$ has a local model
$H_{A'}^1 = \ker d_{A'}^+ / \im d_{A'}$ at $A'$.
In fact the moduli space has the structure of a finite-dimensional
manifold at $A'$, when $d_{A'}$ is injective and $d_{A'}^+$ is surjective.
In that case $H^1_{A'}$ is canonically isomorphic to the tangent space
of the moduli space at $A'$.  These properties are well known for closed
four--manifolds.  A parallel argument also works for  the open case.
Once one fixes a `base' $L^2$ ASD connection $A_0 \in {\mathfrak{M}}$, 
then transversality also works for this case.  
Thus if one takes a generic metric $g'$ with respect to $A_0$, then the
corresponding moduli space $\widetilde{\mathfrak{M}}(A_0, g',w)/{\mathfrak{G}}$
will have the structure of a finite-dimensional manifold.
Notice that these moduli spaces are parameterized by $A_0 \in {\mathfrak{M}}$, and
$g'$ will depend on $A_0$. This causes some difficulty in perturbation theory.

\subsection{Main results}
Here is the main theorem:
\begin{thm} 
Let $S$ be a smooth oriented open four--manifold constructed 
by attaching to a zero handle Casson handles  of homogeneous trees of
bounded type.  Then there exists a complete Riemannian metric of bounded
geometry $g$ and an weight function $w$  on $S$ such that one can
construct  ASD moduli spaces over $S$ as  finite dimensional smooth
manifolds  with respect to $(S,g,w)$.
\end{thm}
This follows from the next two propositions.
In this paper we will introduce admissibility for a pair of a Riemannian
metric and a weight function over an open four--manifold (1.A). Then
we will show the following:
\begin{prop}
Let $S$ be an open four--manifold. Suppose one can equip an admissible pair
$(g, w)$ over $S$. Then one can construct  ASD moduli spaces over
$S$ as  finite-dimensional smooth manifolds.
\end{prop}
In order to obtain admissible pairs over Casson handles of bounded type,
we will use an iterative method.
The main construction in this paper is the following:
\begin{prop} Let $S$ be in  theorem $0.1$.
There exist a  complete Riemannian metric and an  weight function on $S$
so that the  corresponding Atiyah--Hitchin--Singer (AHS)  complex
over $S$ becomes Fredholm. Moreover the metric and the weight function
over  $S$ give an admissible pair.
\end{prop}

Let us  outline  the construction of Riemannian metrics on Casson handles
by restricting to simple cases.  Each  building block of  a Casson handle,
namely a  kinky handle, is diffeomorphic to $\natural (S^1 \times D^3)$
with  two attaching regions; one is a tubular neighborhood of band sums
of Whitehead links (this is connected with the previous block), and the
other  is a disjoint union of the standard open subsets $S^1 \times D^2$
in $\sharp S^1 \times S^2 = \partial (\natural S^1 \times D^3)$ (this
is connected with the next block). The number of end-connected sums
is exactly the number of self-intersections of the immersed two handle.
The simplest Casson handles  have $S^1 \times D^3$ as their building blocks.
We attach a Casson handle to the zero--handle along the attaching circle
and denote it by $S = D^4 \cup CH$.

Let us consider a simple Casson handle, say $CH(\mathbb{R}_+)$, a
periodic Casson handle by positive kinks.  Unlike the Taubes construction,
the building blocks here are open.  In order to make end-connected sums of
building blocks isometrically, one explicitly equips the metrics on the
building blocks. Then the building block as an open manifold becomes an
`open cylindrical' manifold.  As in the Taubes construction, one connects
two attaching regions in a block.  The result becomes a cylindrical
manifold  on which analysis is already well known.  By equipping it with a
suitable weight function, one will apply the  Fourier--Laplace transform
between the cylindrical manifold and its periodic cover.   By a kind of
excision, one obtains a Fredholm AHS complex over  $S$.

This method shows that once one obtains some suitable function spaces
on any open manifolds, then the Fourier--Laplace transform works on their
periodic covers.  We will use this observation iteratively.  In general
a Casson handle can be expressed by an infinite tree with one end point
and with a sign $\pm$ on each edge.
The next simplest Casson handle will be represented as follows.
Let $\mathbb{R}_+$ be the half-line with the vertices $\{0,1,2,\dots\}$.
We prepare another family of  half-lines $\{\mathbb{R}_+^i\}_{i=1,2,\dots}$
assigned with indices.  Then we obtain another infinite tree:
$$R(2) = \mathbb{R}_+ \cup_{i= 1,2, \dots } \mathbb{R}_+^i$$
where we connect $i$ in $\mathbb{R}_+$ with $0$ in $\mathbb{R}_+^i$. 
For example one may assign $-$ on $\mathbb{R}_+$ and $+$ on all $\{
\mathbb{R}_+^i \}_i$.   Then one obtains the corresponding Casson handle
$CH(R(2))$.
In this case the building blocks are diffeomorphic to 
$\overline{\Sigma}_2 = (S^1 \times D^3) \natural (S^1 \times D^3)$
along $\mathbb{R}_+$.
$\overline{\Sigma}_2$  has three attaching components.
One is $\mu$, the tubular neighborhood of the band sum of two Whitehead
links as before.
We will denote the others by $\mu'$ and $\gamma$, where these represent
a generator of $\pi_1(\overline{\Sigma}_2)$.  In order to apply
Fourier--Laplace transform, one  takes  end-connected sums twice.
Firstly  one takes the end-connected sum between $\mu$ and $\mu'$ as before. 
The result is an `open cylindrical' manifold, since  there still
remains one attaching region, $\gamma$. One takes the end-connected
sum of this with $CH(\mathbb{R}_+)$ along $\gamma$.  In this manner,
one obtains  another open manifold, $(\overline{\Sigma}_2 / \mu \sim \mu')
\natural_{\gamma} CH(\mathbb{R}_+)$.  Now we have already two kinds of
analytic preparations.  One is analysis for cylindrical manifolds, and
the other   for the half-periodic Casson handle as we have explained.
By a kind of excision argument, one can verify that the AHS complex on
the open manifold  is Fredholm.  Half of part of its periodic covering is
exactly $CH(R(2))$.  Again by a Fourier--Laplace transform and excision,
one obtains the Fredholm AHS complex over  $S= D^4 \cup CH(R(2))$.

These are simple examples, but the idea works for  much more general
cases of Casson handles.  One may iterate this construction inductively
to more complicated Casson handles.

In the case of end-connected sums, the excision argument tells
us that just the differentials have closed range. In order to see
the finite-dimensionality of the cohomology, we will make explicit
calculations. This is one point where we use de Rham cohomology calculations.
On this point a parallel argument by Seiberg--Witten theory seems to have
some technical difficulty.

\subsection{Directions for further research}
Finally let us indicate some possible developments arising from  this
kind of analysis, assuming perturbation theory.  We would like to propose
here some problems on the study of smooth structures on Casson handles.

Let us consider algebraic surfaces, say the K3 surface.
This decomposes topologically as $2|-E_8| \sharp 3(S^2 \times S^2)$,
and  contains six Casson handles.
One may  guess  that  it would be impossible to do Yang--Mills gauge
theory  on Casson handles inside the K3 surface, or more generally inside
many of algebraic surfaces.
One can verify that at least one of Fredholm theory or perturbation theory
breaks down \cite{kato}.
It seems reasonable to  think that the smooth structure on the Casson
handles in K3 will be so complicated that one might not be able to do
Fredholm theory.  If perturbation theory could work on these, then one
will be able to tell that homogeneous Casson handles inside K3 should
grow more than exponentially.  The argument is  to construct  Yang--Mills
theory on Casson handles and it will lead to a contradiction by
dimension-counting on the moduli spaces.

Let us take two Casson handles, $CH(T_1)$ and $CH(T_2)$
where $T_i$ are the corresponding signed trees.
When $T_1$ is embedded into $T_2$, then there is a smooth embedding,
$CH(T_2) \hookrightarrow CH(T_1)$ preserving the attaching circles,
but one cannot say about converse embeddings of  $CH(T_1)$ into $CH(T_2)$
in general.  The above argument suggests that any CH of bounded type
will not be able to embed into Casson handles  inside K3  preserving
the attaching circles.

In our situation here, one treats Casson handles whose trees grow
polynomially. In fact in our method, one might expect that as Casson
handles grow near exponentially, the continuous spectrum with respect
to the AHS complex will approach zero.  It might be possible that even
exponential growth is already too complicated to obtain Fredholm
theory.  On the other hand one does not know concretely how the signed
trees grow for the case of Casson handles in K3 (for this direction,
see \cite{bizaca}).  In reality, any Casson handles of bounded type
can appear in $S^2 \times S^2 \backslash pt$, and so it would be
interesting to study smooth types on $S^2 \times S^2 \backslash pt$
arising from Casson handles of bounded type.

Next we will consider another problem.  The re-imbedding theorem gives
another Casson handle inside a six stage tower preserving the attaching
circles.  Now any Casson handle $CH$ of bounded type can be embedded into
one of the simplest Casson handles $CH_{\pm} \equiv CH(\mathbb{R}_{\pm})$,
say into $CH_+$.  Let $CH_+(n)$ be the $n$th stage of $CH_+$.  One may
consider a question whether $CH$ can be smoothly embedded into $CH_+(n)$
for some large $n$ preserving the attaching circles.  For this,  we would
like to outline  a possible argument  (see \cite{bizaca and gompf}).  Let $(Z, \partial_0
Z , \partial_1 Z)$ be a smoothly non-product $h$--cobordism with different
Donaldson's polynomials on two boundaries.  By using Kirby calculus
technique, one may find  a decomposition $Z = W \cup U$ where $W$ is
smoothly  product and $CH_+(n)$ appears in both the ends, $\partial_0 Z
\cap W$ and $\partial_1 Z \cap W$.  If $CH$ could be smoothly embedded
into $CH_+(n)$, then one will obtain Yang--Mills moduli spaces over both
of $\partial_i Z \cap W$ whose ends are consisted by CH of bounded type.
Since these are diffeomorphic, the Donaldson's invariants over them  will
have the same numerical value, which would contradict the assumption.
In the above situation, one would be able to conclude that there are no
smooth embeddings of  $CH$ into any finite-stage approximations  $CH_+(n)$
preserving the attaching circles.

The author would like to express his gratitude to the Max Planck
Institut for their hospitality, during his visit.  He also would like
to express his appreciation of the referees for numerous comments and
suggestions.

\section{ Yang--Mills moduli space}
\subsection{Admissible four--manifolds}
In this section, one will construct Yang--Mills moduli spaces over
non compact four--manifolds.  We treat open four--manifolds which can
admit Riemannian metrics and weight functions such that they are able
to construct a Fredholm complex.

Let $Y$ be a non compact smooth four--manifold.  Throughout this section,
one always assumes that $Y$ is simply-connected, and simply-connected
at infinity.

Let $g$ be a complete Riemannian manifold of bounded geometry and 
$w\co Y \to [0, \infty)$ be  a smooth function.
Recall that a complete Riemannian metric is of bounded geometry, if 
(1)\qua the injectivity radius is more than a positive constant $\epsilon >0$ at any 
point,  and   
(2)\qua for any $l \geq 0$, sup$_{x \in Y} |\nabla^l R| < \infty$, where
$\nabla$ and $R$ are with respect to $g$. 
Throughout this paper, we always assume that $g$ is 
of bounded geometry, and $|\nabla^l (w)|C^0(Y) < \infty$ for all $l \geq 1$.

Let $\Lambda^i(Y)$ and $\Lambda^2_+(Y)$ be
the exterior differentials on $i$--forms and self dual $2$--forms with respect to $g$, $i=0,1$. 
Then we have the Atiyah--Hitchin--Singer complex (without coefficient)
as:
$$0 \longrightarrow C_c^{\infty}(Y ;\Lambda^0) \stackrel{d}{\longrightarrow}
  C_c^{\infty}(Y;  \Lambda^1) \stackrel{d^+}{\longrightarrow}
  C_c^{\infty}( Y;  \Lambda^2_+) \longrightarrow 0$$
where $d^+= (1+*) \circ d$.
Using the $L^2$ adjoint operator (here we do not use any weight
functions), we get the next elliptic operator:
$$P= d^+ \oplus d^*\co  C_c^{\infty}( Y; \Lambda^1)
\to C_c^{\infty}( Y;   \Lambda^0 \oplus \Lambda^2_+).$$

 Let us introduce
 weighted Sobolev $k$--norms on $Y$ by:
$$|u|W^k_w = (\Sigma_{l \leq k} \int \exp(w)|\nabla^l u|^2 vol )^{\frac{1}{2}}.$$
 We denote these Sobolev weighted spaces by $W^k_w$ or $L^2_w$ (for $k=0$).
Then  one has a complex of bounded maps:
$$0 \longrightarrow W^{k+2}_w(Y; \Lambda^0) \stackrel{d}{\longrightarrow}
  W^{k+1}_w(Y; \Lambda^1) \stackrel{d^+}{\longrightarrow}
  W^k_w(Y; \Lambda^2_+) \longrightarrow 0$$
Now let us say that the pair $(g,w)$ is \definition{admissible}, if
the following conditions are satisfied:

(1)\qua The above 
consists of a Fredholm complex, namely each differential has closed range, 
and its cohomology group is of finite dimension as a vector space over $\mathbb{R}$.
We denote their cohomology groups by $H^*(\AHS)$, $*=0,1,2$.
Notice that by the condition, one has $L^2_w$ adjoint operators $d^*_w$
and $(d_+)^*_w$.

(2)\qua $Y$ admits a partition $Y = K \cup Y_1 \cup \dots Y_l$ such that
(a) $K$ is a compact subset and each $Y_j$ is an  open subset of $Y$, and
(b) let us put $D = d^*_w \oplus d_+$ or 
   $d \oplus (d_+)^*_w$.
Then there is a positive constant $C_k$ such that 
for any $u \in C_c^{\infty}(Y_j; \Lambda^*)$, $*=1,2$, 
one has the bound $|u|W^{k+1}_w \leq C|D(u)|W^k_w$.

(3)\qua For any $f \in C^1(Y)$ with $|df|L^2_w(Y) < \infty$,
there is $\overline{f} \in \mathbb{R}$ with $f- \overline{f} \in L^2_w(Y)$.

(4)\qua There is a compact subset $K \subset Y$ and a homeomorphism 
$Y \backslash K \cong S^3 \times [0, \infty)$.

In sections $2$ to $5$, we will construct admissible pairs
for all Casson handles constructed from
homogeneous signed trees  of bounded type.
\begin{prop}  Let $(g,w)$ be an admissible pair on $Y$. 
Then    $H^0(\AHS)=H^1(\AHS)=0$ and
$H^2(\AHS) \geq k = H^2_+(Y; \mathbb{R})$.
\end{prop}
\begin{proof}
Clearly  $H^0(\AHS)=H^1(\AHS)=0$ by the admissibility condition $(3)$.
We show $\dim H^2 (\AHS)\geq k$. For this, one takes two steps.

Let $H^2_{\text{cp}}(Y;\mathbb{R})^+$ be a subspace of 
$H^2_{\text{cp}}(Y;\mathbb{R})$ consisting of
vectors with  $ \int_Y  u \wedge u  >0$.
Then  $H^2_{\text{cp}}(Y;\mathbb{R})^+$ is 
a linear subspace of $H^2(Y;\mathbb{R})$ of 
dimension $k$, since the natural map
$I\co H^2_{\text{cp}}(Y; \mathbb{R}) \to H^2(Y ;\mathbb{R})$
 gives an isomorphism.
 
Now
by the above, one has $\dim H^2_{\text{cp}}(Y;\mathbb{R})^+ =k$.
For any element $u \in C^{\infty}(Y; \Lambda^2)$, 
let us denote by $u^+$ the projection to $\Lambda^2_+$ part.
Then one defines:
$$p\co H^2_{\text{cp}}(Y;\mathbb{R})^+ \to H^2(\AHS)$$
by assigning $[u] \mapsto [u^+]$.
This map is well defined. We show that $p$ is an injection.

Suppose $[u^+] =0 \in H^2(\AHS)$. Then one has $\alpha \in W^1_w(Y; \Lambda^1 )$
such that $u^+ =d^+(\alpha)$, or $(u-d(\alpha))^+ =0$.
Let us put $v= u-d(\alpha)$, and take the cup product $\langle v,v \rangle = \int_Y v \wedge v$.
One can use  the Stokes theorem to see
$\int_Y u \wedge d(\alpha)  = \int_Y d( u \wedge \alpha)    =0$  ($u$ is a closed form).
Thus one has the equality $ \int_Y v \wedge v = \int_Y u \wedge u \geq 0$.
On the other hand since $v$ does not have self-dual part, one has $ \int_Y v \wedge v \leq 0$.
This shows $u = d(\alpha)$. This implies $[u]=0 \in H^2_{\text{cp}}(Y;\mathbb{R})^+$ since
$I$ gives an isomorphism. 
Thus one concludes $\dim H^2  \geq k$.
\end{proof}

\begin{exmp}
{\rm In general, $\dim H^2_{cp}(Y;\mathbb{R})^+$
and $\dim H^2(\AHS) = \ker (d_+)^*_w$  do not coincide. One takes
$Y = \mathbb{R}^2 \times \mathbb{R}^2$, where one equips the
standard metric on $\mathbb{R}^2$, and the product one on $Y$.
Let us choose a weight function $w$ on $Y$.
We show $\dim H^2(\AHS)$ is nonzero.
Let $u_1$ and $u_2$ be $2$--forms on $\mathbb{R}^2$
respectively, such that $|u_i| \equiv 1$, $i=1,2$,
pointwisely. Then clearly $du_i=0, i=1,2$, and
$u_1+u_2$ gives a self-dual $2$--form
with bounded pointwise norm. 
Then one puts $v= \exp(-w)(u_1+u_2) \in L^2_w(Y; \Lambda^2_+)$.
Moreover $v$ satisfies the equation $d(\exp(w)v) =0$. 
Thus $v$ is a nontrivial element in $H^2(\AHS)$.}
\end{exmp}

\subsection{ASD moduli space}
Let $E \to Y$ be a $G$--vector bundle (where $G=SO(3)$ or $U(1)$) such
that except
on some compact subset $K \subset Y$, $E|Y \backslash K$
is trivial. One denotes the corresponding principal $G$--bundle by $P$.
For $SO(3)$--bundles, $E$ is determined by $w_2(E) \in H^2(Y;\mathbb{Z}_2)$.
In later sections, we fix  a trivialization of $E| Y \backslash K$.
Thus we fix an $SO(3)$--bundle with $w_2$ and $p_1(E)$.
Let $A$ be a smooth connection over $E$ such that except some compact subset on $Y$,
it satisfies the ASD equation:
$$F_A +*F_A=0$$
where $F_A$ is the curvature form of $A$. 
Let us denote by $\widetilde{R}(Y)$ the set of smooth connections as above satisfying 
$|F_A|L^2(Y) < \infty$. By changing a trivialization, one may assume
$|a|L^2(Y \backslash K) < \infty$ where $A = d+a$ on $Y \backslash K$.
Later we always assume this property.
If $A$ is the trivial connection except some compact subset, then it is an
element in  $\widetilde{R}(Y)$.
\begin{lem}
For $A \in \widetilde{R}(Y)$,  $p_1(A) = \frac{1}{4 \pi^2} \int_Y tr(F_A \wedge F_A)$ is an integer.
\end{lem}
In order to show this, one uses the following:
\begin{sublem}{\rm\cite{freed and uhlenbeck}}\qua
Let $D \subset Y$ be a $\delta$--ball with any point as the center, where
$\delta >0$ is sufficiently small.
 Then there exists another small $\epsilon >0$ with the following property;
suppose $|F_A|L^2(D) < \epsilon$. Then there exists a gauge transformation $g$
over $D$, such that: 
$$\textstyle{\sup_{x \in D}} | \nabla^l( g^*(A) -d)|(x) \leq  C_l |F_A|L^2(D)$$
where $C_l$ are constants, independent of $A$.
\end{sublem}

{\bf Proof of lemma}\qua
Let us take any small $\epsilon >0$. Then there exists a compact subset
$K \subset Y$, and a smooth gauge transformation
$g \in C^{\infty}(Y; \Aut(E))$ such that
$|g^*A -d|W^1(Y \backslash K) < \epsilon$, 
where $d$ is the trivial connection ($Y$ is simply-connected at infinity).
In order to get $g$, one needs to use the above sublemma and the patching
argument.  We omit to describe the process (see 
\cite{dk}).
One may assume that $K$ is a compact submanifold with boundary.
Let $D(K)$ be a double of $K$. Then using a cut-off function 
$\varphi$ around $\partial K$,
one may extend  $A'  \equiv \varphi(g^*A )$ over $D(K)$, by putting $d$ 
over the other side of $K \subset D(K)$.
By the Chern--Weil theory, $p_1(A')$ is an integer.
On the other hand, since $p_1(A) - p_1(A')$ may be arbitrarily small as above,
one concludes that $p_1(A)$ is an integer.
This completes the proof.
\endproof

Let us put $\Ad(P) = P \times_G \mathfrak{G}$, where 
$\mathfrak{G}$ is the Lie algebra of $G$.  ($P \times_G \mathbb{R}^3 =E$).
We also denote by $\mathfrak{G}$ the trivial Lie $G$--bundle.
Then we have the Atiyah--Hitchin--Singer complex (AHS complex) as:
$$0  \longrightarrow C_c^{\infty}(Y; \Ad(P))  \stackrel{d_A}{\longrightarrow}  
 C_c^{\infty}(Y; \Ad(P)\otimes \Lambda^1) \stackrel{d^+_A}{\longrightarrow}
 C_c^{\infty}(Y; \Ad(P)\otimes\Lambda^2_+) \longrightarrow 0$$
where $d^+_A= (1+*) \circ d_A$.
Using the $L^2_w$ adjoint operator, we get the next elliptic operator:
$$P(A)_w = ( d_A)_w^* \oplus  d_A^+\co W^{k+1}_w(Y; \Ad(P) \otimes \Lambda^1)
\to  W^k_w(Y; \Ad(P) \otimes ( \Lambda^0 \oplus \Lambda^2_+)).$$

Let us fix $A_0 \in \widetilde{R}(Y)$.  One defines an affine Hilbert space as:
$${\mathfrak{A}}_k (P)_0=  \{A_0 +a| a \in W^k_w(Y; \Ad P \otimes T^*Y) \},
  \qua k \geq 3.$$ 
Let us take $g \in  C^1_{\text{loc}}(Y;\Aut E)$. By embedding as 
$g \in 
C^1_{\text{loc}}(Y; \Hom(E,E))$, one may consider 
$\nabla_{A_0}g \in C^0(Y;\Hom(E,E)\otimes T^*Y)$.
Notice that if $g$ is locally $W^4$, then it is   of $C^1$ class.
Now one defines the weighted Sobolev gauge group:
\begin{equation}\begin{aligned}
& {\mathfrak{G}}_l(P)= \{ h \in W^l_{loc}(Y; \Aut(E))|\nabla_{A_0}h \in
W^{l-1}_w(Y; \Hom(E,E) \otimes T^*Y) \}, \qua l \geq 4. \\
&{\mathfrak{G}}_l(P)_0= \{ h \in W^l_{loc}(Y; \Aut(E))|
h - id \in W^{l-1}_w(Y; \Hom(E,E) \otimes T^*Y) \}.
\end{aligned}\end{equation}
First of all one has the following property of $\mathfrak{G}_l(P)$.
Suppose one of the following two conditions;
(1) $G=U(1)$, or (2) $G= SO(3)$ and $A|Y \backslash K 
\in W^k_w(Y \backslash K ; \Ad(P) \otimes \Lambda^1)$
where  the trivialization of $E|Y \backslash K$ is fixed. 
Then:
\begin{lem}
For each $h \in {\mathfrak{G}}_{k+1}(P)$, there exists $\overline{h} \in G$ such that
$h-\overline{h} \in W^k_w(Y \backslash K; \Aut(E))$, where $\overline{h}$ is a constant gauge
transformation.
\end{lem}
\begin{proof}
This follows from Kato's inequality, 
$|\nabla_{A_0}h| \geq d|h|$
almost everywhere, and admissibility condition (3).
\end{proof}

The Lie algebras of ${\mathfrak{G}}_l(P)$ and ${\mathfrak{G}}_l(P)_0$ are 
correspondingly as follows:
\begin{equation}\begin{aligned}
& {\mathfrak{g}}_l(P)= \{ h \in W^l_{loc}(Y; \Ad(P))|\nabla_{A_0}h \in
W^{l-1}_w(Y; \Ad(P) \otimes T^*Y) \},\\
& {\mathfrak{g}}_l(P)_0 = W^l_w(Y; \Ad(P)).
\end{aligned}\end{equation}
\begin{lem} {\rm(1)}\qua  ${\mathfrak{G}}_{k+1}(P)_0$ acts 
on ${\mathfrak{A}}_k(P)_0$ by
$g^*(A_0+a) = g^{-1} \nabla_{A_0} g + g^{-1} ag$.

{\rm(2)}\qua Suppose $P$ is a $U(1)$--bundle with $G= U(1)$.
Then ${\mathfrak{G}}_{k+1}(P)$ acts on
${\mathfrak{A}}_k(P)$.
\end{lem}
The proof is standard, and we omit it.
\vspace{2mm}

Notice that  ${\mathfrak{G}}_{k+1}(P)$ may not act on ${\mathfrak{A}}_k(P)_0$,
since $A_0$ may live only in $L^2(Y)$.

Let us define:
$$\widehat{{\mathfrak{M}}}_k(P) =
  \{ A \in {\mathfrak{A}}_k(P)_0: F_A+ *F_A =0\}.$$

Then ${\mathfrak{G}}_{k+1}(P)_0$ acts on 
$\widehat{\mathfrak{M}}_k(P)$.  Thus one gets the quotient space:
$${\mathfrak{M}}_k(P) = \widehat{{\mathfrak{M}}}_k(P)/ {\mathfrak{G}}_{k+1}(P)_0.$$
\begin{rem}
\textrm{
${\mathfrak{M}}_k(P)$ is constructed after choosing 
a base connection $A_0 \in \widetilde{R}(Y)$ and $W^k_w$ Hilbert spaces.
It can be shown that this space is in fact independent of choice
of $k$. However it would definitely depend on choice of
$A_0$ ($A_0$ lies only in $L^2$ with respect to the trivialization
on the end). Thus one could express the space as:
$${\mathfrak{M}}_k(P) = {\mathfrak{M}}_k(P, A_0).$$}
\end{rem}
 
Let us find a linear space which corresponds to the tangent space
of ${\mathfrak{M}}(P)_k$.
\begin{lem}
There exist natural continuous maps:
$$\exp\co   {\mathfrak{g}}_{k+1}(P) , \quad
  ( W^{k+1}_w(Y; \Ad (P)))  \to {\mathfrak{G}}_{k+1}(P) \quad
 ({\mathfrak{G}}_{k+1}(P)_0).$$
Moreover, let us put 
${\mathfrak{G}}'_{k+1}(P)(\epsilon) =
\{ g \in {\mathfrak{G}}_{k+1}(P)| |g-1 | W^{k+1}_w(Y)  \leq \epsilon \}$.
Then for sufficiently small $\epsilon >0$,
there exists $\log\co {\mathfrak{G}}_{k+1}(P)_0(\epsilon) \to W^{k+1}_w(Y;
\Ad P)$
which inverts $\exp$.
One has a similar statement for 
$\exp| W^{k+1}_w(Y; \Ad (P))$.
\end{lem}
We also omit the proof.

Let us take any $A \in {\mathfrak{A}}_k(P)_0$, and
 consider the following continuous map:
$$d_A\co {\mathfrak{g}}_{k+1}(P)  \to W^k_w(Y; \Ad P \otimes T^*Y)$$
by $u \to \frac{d}{dt}(\exp(tu)A)|_{t=0}$. Locally this map is expressed
as $u \to du+[u,a]$. 
\begin{lem} Suppose an admissible pair $(g, w)$ is given. Then
for the AHS complex, $d_A$ and $d_A^+$  are  closed operators.
\end{lem}
\begin{proof}
Let us show that $d_A| W^{k+1}_w(Y; \Ad(P))$ has closed range.
Let $A_0' \in {\mathfrak{A}}_k(P)$ be a smooth connection which coincides
with $A$ on
some compact subset $K$, and  is trivial except on another compact subset.
Then one has $A=A_0'+a$, $a \in W^k(Y; \Ad P \otimes T^*Y)$.
Let $\varphi$ be a cut-off function such that for some compact subsets
$K \subset K' \subset K'' \subset  Y$, one has $\varphi|K' \equiv 1$
and $\varphi|(K'')^c  \equiv 0$.  By the assumption, $d\co W^{k+1}_w(Y;
\Lambda^0) \to W^k_w(Y; \Lambda^1)$ is invertible (satisfies $|d(u)|W^k_w
\geq C|u|W^{k+1}_w$).  Then one sees $d_{(1-\varphi) a}$ is invertible
by choosing sufficiently large $K' \subset K''$, using the Sobolev
embedding, $W^k_{\text{loc}} \hookrightarrow C^0$.
Then one has the following estimate:
\begin{equation}\begin{aligned}
|u|W^{k+1}_w & \leq C( |d_A(u)|W^k_w + |u|W^k_w) \\
& \leq C( |d_A(u)|W^k_w + |\varphi u|W^k_w +|(1-\varphi) u|W^k_w) \\
 & \leq C(|d_{A}(u)|W^k_w  + |\varphi u|W^k_w + |d_{(1- \varphi) a }((1-\varphi )u)|W^k_w) \\
& \leq  C(|d_{A}(u)|W^k_w  + |\varphi u|W^k_w + 
|(1-\varphi )d_{(1-\varphi) a}(u)|W^k_w + \\
 & \qquad |[d_{(1- \varphi) a}, (1-\varphi)]u|W^k_w)  \\
&  \leq C(|d_{A}(u)|W^k_w  + |\varphi u|W^k_w + 
|(1-\varphi )d_A(u)|W^k_w + \\
& \qquad |[d_{(1- \varphi) a}, (1-\varphi)]u|W^k_w +  
|(1-\varphi )[\varphi a + A_0', u]|W^k_w) \\
& \leq   C(|d_{A}(u)|W^k_w  + |\psi u|W^k_w) 
 \end{aligned}\end{equation}
where $\psi $ is another compactly-supported cut-off function with 
$\psi|\Supp\varphi \equiv 1$.
The last inequality shows that $d_A$ is a closed operator (see lemma $3.1$).

Next let us consider $d_A^+$.  Notice that we have now the adjoint operator $(d_A)^*_w$.
Then it is enough to see $d_A^+( \ker (d_A)^*_w \cap W^{k+1}_w) \subset
W^k_w(Y; \Ad(P) \otimes \Lambda^2_+)$
has closed range (see lemma $1.6$ and lemma $3.1$).
Let us put $D= (d_A)^*_w \oplus d_A^+$. Then using admissibility condition (2)
and the above argument, one gets a similar inequality
$|u|W^{k+1}_w \leq C( |D(u)|W^k_w + |\psi u|W^k_w)$
where a cut-off function $\psi$ has compact support.
In particular  $D\co W^{k+1}_w \to W^k_w$ has closed range.

Let us take $u \in W^{k+1}_w(Y; \Ad(P) \otimes \Lambda^1)$ with
$(d_A)^*_w(u) =0$.  Then by the above two facts,
one has a similar  inequality:
$$|u|W^{k+1}_w \leq C( |d_A^+(u)|W^k_w + |\psi u|W^k_w).$$
This shows that $d_A^+$ also has closed range.
This completes the proof.
\end{proof}
\begin{cor}
Let $A \in {\mathfrak{A}}_k(P)_0$. Then the AHS complex
\begin{multline*}
0 \longrightarrow  W^{k+1}_w(Y; \Ad(P)) \stackrel{ d_A }{\longrightarrow}  
  W^k_w(Y; \Ad(P) \otimes \Lambda^1) \\
\stackrel{ d_A^+ }{\longrightarrow} W^{k-1}_w(Y; \Ad(P) \otimes \Lambda^2_+)
  \longrightarrow  0
\end{multline*}
is a Fredholm complex of index $-2p_1(P) + 3 \left( \dim H^1(\AHS)
- \dim H^2(\AHS)\right)$.
\end{cor}
The index computation uses the excision principle, or relative index
theorem \cite{gromov and lawson}.  One can also use \cite{atiyah patodi
and singer} and the method
in section 5 here.  One denotes their cohomology groups by $H^*_A(\AHS)$,
$*=0,1,2$.

Let us consider the restriction:
$$d_A\co W^{k+1}_w(Y; \Ad P ) \to W^k_w(Y; \Ad P \otimes T^*Y).$$
By the above corollary,  one may consider the adjoint operator:
$$(d_A)^*_w\co W^k_w(Y; \Ad P \otimes T^*Y) \to W^{k-1}_w(Y; \Ad P).$$
One understands  $(d_A)^*_w$ in the geometric way as follows.
Let us  take  $A \in {\mathfrak{A}}_k(P)_0$.
Then one has the `tangent space' of 
$\mathfrak{A}_k(P)_0 / \mathfrak{G}_{k+1}(P)_0$ at $[A]$:
$$\frac{W^k_w(Y; \Ad P \otimes T^*Y)}{d_A(W_w^{k+1}(Y; \Ad P))}\cong 
\ker(d_A)^*_w\co W^k_w(Y; \Ad P \otimes T^*Y) \to W^{k-1}_w(Y; \Ad P).$$

Let us take a smooth family of connections $A(t) \in {\mathfrak{A}}_k(P)_0$,
$t \in [0,1]$. Then  one has:
$$\left.\frac{\partial F_{A(t)}^+}{\partial t }\right|_{t=0}
  = (d_{A(0)} + *d_{A(0)})B \equiv d_{A(0)}^+B$$
where $B=\left.\frac{d A(t)}{dt}\right|_{t=0} \in W^k_w(Y; \Ad P \otimes T^*Y)$.
Notice that near infinity, $F_{A(t)}$ can be written as $dA(t)+ A(t) \wedge A(t)$.
Each $A(t)$ can be written as $A(t) = A_0 + a_t$, where 
$a_t \in W^k_w(Y; \Ad P \otimes \Lambda^1)$.
Since $F_{A_0}^+ =0$, 
one sees that $F_{A_t}^+$ is a smooth family in $W^{k-1}_w(Y; \Ad P \otimes \Lambda^2_+)$.
Now one gets the `tangent space' of ${\mathfrak{M}}_k(P) $ at $[A]$ as:
$$\ker~(d_A)^*_w \oplus d_A^+\co W^k_w(Y; \Ad P \otimes T^*Y) \to
W^{k-1}_w(Y; \Ad P \otimes (\Lambda^0 \oplus \Lambda^2_+)).$$
In particular this space is isomorphic to $H^1_A(\AHS)$ which is of finite dimension.
\begin{lem}
$\ker d_A\co W^{k+1}_w(Y; \Ad(P)) \to W^k_w(Y; \Ad(P) \otimes \Lambda^1)$ is zero.
\end{lem}
\begin{proof}
This follows from Kato's inequality, $|\nabla_A u| \geq |d|u||$ almost
everywhere.
Notice that there are also no $L^2$ functions $u$ with $d_A(u)=0$.
This is also seen as follows: Suppose $u \in L^2(Y; \Ad(P))$, $|u|L^2=1$, 
satisfies $d_A(u)=0$. Then for  $g= \exp(u) \in \mathfrak{G}$,
one has $g^*(A)=A$. Then it is known that  for any path $l$ between $p$ and $q$
 in $Y$, one has $P_l \circ g(p) =g(q) \circ P_l$, where $P_l$ is the parallel
translation along $l$.
Let $ q \to \infty$. Then $g(q)$ approaches to the identity. On the other hand
since $P_l$ holds the inner product, one has $\langle P_l(g(p)(u)),
g(q)(P_l(u)) \rangle_q \sim \langle g(p)(u), u \rangle_p$. This shows
that $g(p)$ is near the identity. This impossible.
This completes the proof.
\end{proof}

Let us consider
 ${\mathfrak{M}}_k(P) = \widehat{{\mathfrak{M}}}_k(P)_0 / {\mathfrak{G}}_{k+1}(P)_0$, 
and show that under some assumptions, this 
space is a finite-dimensional smooth manifold.
For $A \in {\mathfrak{A}}_k(P)_0$, let us put:
$$CG_A^k=  \{A' \in {\mathfrak{A}}_k(P)_0|(d_A)^*_w(A' -A)=0\}. $$
$$I\co CG_A^k \times {\mathfrak{G}}_{k+1}(P)_0 \to 
     {\mathfrak{A}}_k(P)_0;~~ (A, g) \mapsto  g^*A.$$
Let us calculate $dI$ at $(A,\id)$. If $u \in W^{k+1}_w(Y; \Ad(P))$ and
$v \in W^k_w(Y; \Ad(P) \otimes \Lambda^1)$, one has:
$$\frac{d}{dt}(\exp  tu)^*(A+tv)|_{t=0} = d_Au+v.$$
Thus $dI$ is given by:
\begin{equation}\begin{aligned}
dI=1 \oplus d_A\co
\ker(d_A)^*_w \cap & W^k_w(Y; \Ad P \otimes T^*Y) \times 
W_w^{k+1}(Y; \Ad P ) \\
 & \to
W^k_w(Y; \Ad P \otimes T^*Y).
\end{aligned}\end{equation}
Since $d_A$ has closed range, it is clear that the above map is an isomorphism  at every $(A,\id)$. 
\begin{cor}
There are neighborhoods $U \subset {\mathfrak{A}}_k(P)_0$ of $A$
and $V \subset {\mathfrak{G}}_{k+1}(P)_0$ of $\id$ such that:
$$I\co U \cap CG_A^k \times V \to {\mathfrak{A}}_k(P)_0$$
is a homeomorphism into its image.
By restriction,
$$I\co U \cap CG_A^k(P) \cap \widehat{\mathfrak{M}}_k(P)_0 \times V \to \widehat{\mathfrak{M}}_k(P)_0$$
is also a homeomorphism into $\widehat{\mathfrak{M}}_k(P)_0$.
\end{cor}
Now one has the following:
\begin{prop}
Suppose  $d_A^+\co W^k_w(Y; \Ad P \otimes T^*Y) \to W^{k-1}_w(Y; \Ad P \otimes \Lambda_+^2)$
is a surjection for  $A \in \widehat{\mathfrak{M}}_k(P)_0$.
Then ${\mathfrak{M}}_k(P) = \widehat{\mathfrak{M}}_k(P)_0/ {\mathfrak{G}}_{k+1}(P)_0$ is  a  finite-dimensional smooth manifold
at $[A]$.  Its tangent space is naturally isomorphic to $\ker d_A^+/ \im d_A$.
\end{prop}
\begin{proof} This follows from the general facts on group actions
\cite[page 48]{freed and uhlenbeck}.
From the above corollary, it is enough to verify that the above $I$
gives a slice for the projection:
$$\widehat{\mathfrak{M}}_k(P)_0 \to {\mathfrak{M}}_k(P) =
  \widehat{\mathfrak{M}}_k(P)_0 / {\mathfrak{G}}_{k+1}(P)_0.$$
These can be verified by bootstrapping  as \cite{freed and uhlenbeck}.
We omit to write it.
\end{proof}

\subsection{Perturbation of Riemannian metrics}
In order to construct smooth moduli spaces,
one uses K\,Uhlenbeck's generic metric theorem.
Let $(g,w)$ be an admissible pair for $Y$.
Let us choose a smooth map $h\co Y \to [0, \infty)$ 
with $h(x) \geq  \frac{w}{2}(x)$. Then    one 
introduces the following Banach manifold:
$${\mathfrak{C}} = \{ \phi \in C^l(Gl(TM)):
\lim \textstyle{\sup_K} ( \Sigma_{j=0}^l e^{h} \left| \nabla^j (\phi^*g -
g)\right| K = 0,  K \text{ compact} \}.$$
Let us take $\phi \in \mathfrak{C}$, and put $g'=\phi^*g$.
Then one has the  AHS complex with respect to $g'$:
$$0  \longrightarrow  W^{k+1}_w((Y,g'); \Lambda^0)
  \stackrel{ d }{\longrightarrow} W^k_w((Y,g'); \Lambda^1)
  \stackrel{ d^+ }{\longrightarrow}  W^{k-1}_w((Y,g'); \Lambda^2_+)
  \longrightarrow  0.$$
One sees that this is a Fredholm complex with the same index as the unperturbed one.
Our aim in this section is to show the following:
\begin{prop}
Suppose $Y$ is indefinite. 
Then by a small perturbation of the Riemannian metric,
one has no orbit of reducible connections in ${\mathfrak{M}}_k(P)$.
In particular ${\mathfrak{M}}_k(P)$ is a smooth manifold.
\end{prop}
In order to verify this, we follow \cite{freed and uhlenbeck}.
Let us take any $L^2$ ASD connection $A_0$ over $Y$, and 
consider ${\mathfrak{A}}_k(P)_0 = \{ A_0 +a | a \in  W^k_w(Y; \Ad(P) \otimes \Lambda^1)$.
Then one introduces a ${\mathfrak{G}}_{k+1}(P)_0$ equivariant map:
$$ P_+\co {\mathfrak{A}}_k(P)_0 \times {\mathfrak{C}}
  \to W^{k-1}_w((Y,g); \Lambda^2_+)$$
by $P_+(A, \phi) = P_+(g)(\phi^{-1}(F_A))$,
where $P_+(g)$ is the projection to the self dual part with respect to $g$. 
Notice that this map is well-defined since $A_0$ is $L^2$ ASD with
respect to $g$.  Let us put $\overline{{\mathfrak{M}}}_k(P) = P_+^{-1}(0)$.
\begin{prop}{\rm\cite{freed and uhlenbeck}}\qua
$\overline{{\mathfrak{M}}}_k(P) \cap {\mathfrak{A}}_k^*(P) \times {\mathfrak{C}} /  {\mathfrak{G}}_{k+1}(P)_0$
is a smooth Banach manifold.
\end{prop}
\begin{proof}
We sketch its proof.
First we see that $d P_+$ is surjective at any $(A, \varphi)$ with
$P_+(A, \varphi) =0$.
Then it follows that $\overline{{\mathfrak{M}}}_k(P)$ is a Banach manifold
on which ${\mathfrak{G}}_{k+1}(P)_0$ acts.
Then as before by making a slice for the action, one gets the result.

Let ${\mathfrak{c}}$ be the Lie algebra of ${\mathfrak{C}}$. Then $d P_+$
splits as:
$$d P_+ = d_1 P_+ \oplus d_2 P_+\co
  W^k_w(Y; \Ad(P) \otimes \Lambda^1 ) \oplus {\mathfrak{c}} \to 
  W^{k-1}_w(Y ; \Ad(P) \otimes \Lambda^2_+ ),$$
where $d_1 P_+(\alpha) = P_+(g)(\varphi^{-1}(d_A(\alpha)))|(u, \varphi)$, and
$d_2 P_+(r) = P_+((\varphi^{-1})^*(r^*F))$.
We show that  the differential of $P_+$ is surjective --
notice that
$$\ker d_A \cap W^{k+1}_w(Y; \Ad(P)) = 0.$$

Let us consider the AHS complex:
$$0 \longrightarrow  W^{k+1}_w( \Ad(P) )  \stackrel{ d_A }{\longrightarrow}  
  W^k_w( \Ad(P) \otimes \Lambda^1 )  \stackrel{ d_A^+ }{\longrightarrow}
  W^{k-1}_w(\Ad(P) \otimes \Lambda^2_+)  \longrightarrow  0.$$
Since this is Fredholm, one sees $d P_+$ has finite codimension.
Let us take a representative $u \in \coker d  P_+$ with
 $(d_A^+)^*(u)=0$. Then one has $(d_A^+)^*(e^w u )=0$.
Then one has the equations, $d_A(F_A)=d_A^*(F_A)=
d_A(v) =d_A^*(v) =0$, where $v=\varphi^*(e^w u)$.
Then the same argument as  \cite[page 56]{freed and uhlenbeck} shows that
on an open dense subset of $Y$,
$F_A$ can be expressed as $\alpha \otimes a \in \Lambda^2_+ \otimes \Ad(P)$,
with $|a| =1$ (pointwise norm) and $d_A(a)=0$. 
By the irreducibility, it follows $a=0$, which  contradicts to non
triviality of $p_1(P)$.  This completes the proof.
\end{proof}
\begin{prop}
Let us fix $A_0$, an $L^2$ ASD connection with respect to $g$.
Then for a Baire set of $\phi \in {\mathfrak{C}}$, there are
no reducible connections in ${\mathfrak{A}}_k(P)_0$ with respect to
$\phi^*(g)$.
\end{prop}
\begin{proof}
Suppose $A_0$ is reducible and
denote the corresponding $U(1)$--connection by the same $A_0$.
Let $P$ be a $G = U(1)$--bundle. Then, similarly to before, one puts: 
\begin{equation}\begin{aligned}
& {\mathfrak{A}}_k(P)_0 = \{ A_0 + a| a \in W^k_w(Y; \Lambda^1) \}, \\
& {\mathfrak{G}}_{k+1} (P)_0= \{ h \in W^{k+1}_{loc}(Y; \Aut(E))|
     h - \id \in W^{k+1}_w(Y;  T^*Y \otimes \mathbb{C}) \}.
\end{aligned}\end{equation}
Then ${\mathfrak{G}}_{k+1}(P)_0$ acts on ${\mathfrak{A}}_k(P)_0$.
Now, as above, one considers $P_+\co {\mathfrak{A}}_k(P)_0
\times {\mathfrak{C}} \to W^{k-1}_w((Y, g'); \Lambda^2_+)$.
As \cite{freed and uhlenbeck}, if $d P_+$ is  not surjective, then it
follows $F_{A_0}=0$.
This shows that $dP_+$ is surjective.

Now
one has a smooth Banach manifold
$\overline{{\mathfrak{M}}}_k(P) = P_+^{-1}(0) / {\mathfrak{G}}_{k+1}(P)_0$.
Let us consider a Fredholm map between Banach manifolds
$\overline{\pi}\co \overline{{\mathfrak{M}}}_k(P) \to {\mathfrak{C}}$.
Its Fredholm index is the one of the following:
$$0 \longrightarrow  W^{k+1}_w(Y; \Lambda^0 )
  \stackrel{ d }{\longrightarrow}  W^k_w(Y;  \Lambda^1 )
  \stackrel{ d^+ }{\longrightarrow}  W^{k-1}_w(Y; \Lambda^2_+)
  \longrightarrow  0.$$
One has $H^0(\AHS)=H^1(\AHS)=0$, and $H^2(\AHS) \geq b^2_+$. This is also the case 
when one perturbs the Riemannian metrics slightly.
In particular, if $b^2_+ >0$, then one has $\dim H^1 - \dim H^2 < 0$.
Then the Sard--Smale theorem  shows that for a Baire set of
$\mathfrak{C}$, 
there are no
$L^2$ ASD connections over nontrivial line bundles.
This shows that there are no orbits of reducible ASD connections in
${\mathfrak{M}}_k(P)$, for a Baire set of $\mathfrak{C}$.
This completes the proof.
\end{proof}

Let us put the projection
$\pi\co \overline{\mathfrak{M}}_k(P) \to {\mathfrak{C}}$.
Then the direct application of \cite{freed and uhlenbeck} to this case shows:
\begin{cor}{\rm\cite{freed and uhlenbeck}}\qua Suppose $b^2_+ >0$. Then
a Baire set of $\phi \in {\mathfrak{C}}$ exists such that 
$\overline{ {\mathfrak{M}}}_k(\phi) \equiv \pi^{-1}( \phi)$  are smooth finite dimensional
manifolds.
\end{cor}

\section{A complete Riemannian metric on int $W_0$}
Roughly speaking the set of the anti-self-dual connections modulo
equivalence forms a finite-dimensional manifold, in the case when
the base four--manifold is closed. Then dimension is equal to the Fredholm 
index of the elliptic complex, the Atiyah--Hitchin--Singer  complex.
One may make a Fredholm operator  for a four--manifold with
boundary, using suitable complete Riemannian metrics on its interior.
In this case one uses the weighted Sobolev spaces on the elliptic
operator. For our purpose this is a basic construction, and later
we will use infinitely many such open Riemannian  manifolds.  
In fact by gluing out all of them iteratively, we will get 
another complete Riemannian manifold such that the Atiyah--Hitchin--Singer
complex becomes Fredholm between the weighted Sobolev spaces.
The aim in later sections is to construct admissible pairs
over Casson handles of bounded type.

\subsection{A specified embedding  of $W_0$ into $\mathbb{R}^N$}
Let $W_0$ be a compact manifold with boundary of arbitrary dimension $n$.
Let us choose compact submanifolds with boundary $M_1, N_0 \subset
\partial W_0$  of the same
dimension, with empty intersection, $M_1 \cap N_0 = \emptyset$.
Suppose $M_1$ and $N_0$ are diffeomorphic to each other, and denote
$\widetilde{W}_0 = \interior W_0 \cup M_1 \cup N_0$. 
Let us take a countable number of $\widetilde{W}_0$ and index
them by $\widetilde{W}_0^i = \widetilde{W}_0$,  $i \in \bf{Z}$. Then
one constructs an open manifold  $Y$  by $\cup \widetilde{W}_0^i$
identifying $M_1^i$ with $N_0^{i+1}$:
$$Y = \dots \widetilde{W}_0^i \cup_{M_1^i \sim N_0^{i+1}} \widetilde{W}_0^{i+1} \cup \dots$$
There is a natural smooth $\mathbb{Z}$--action on $Y$, and we call $Y$
a \definition{periodic open manifold}.
There are obvious notions of periodic vector bundles over $Y$, and
periodic differential operators over $E$.

In order to define (weighted) Sobolev spaces over $Y$ we need, first
of all, a complete Riemannian metric on $\interior W_0$ which extends
to the one on $Y$.  We recall an explicit  picture of end-connected sums
along $M_1$ and $N_0$.  Let us choose  smooth embeddings:
$$\gamma (\gamma')\co M_1  (N_0)
\times [0, \infty) \hookrightarrow \partial W_0 \times [0, \infty)
\subset \interior W_0.$$ 
Then we have the following equality (at this stage it is $C^{\infty}$ sense, by smoothing corners):
$$\interior W_0 \natural \interior  W_0 =
\interior W_0 \backslash \im \gamma 
\cup \interior W_0 \backslash \im \gamma '$$
where we identify $\partial M_1 \times [0, \infty) \cup M_1 \times \{0\}$ with
$\partial N_0 \times [0, \infty) \cup N_0 \times \{0\}$.
In order to make this construction compatible with the metric on
$\interior W_0$, we choose a particular style of Riemannian metric on
int $W_0$  as follows.  Let us choose an embedding:
$$I\co W_0 \hookrightarrow \mathbb{R}^{N+2}$$
with the following six properties:

(1)\qua Let us choose a small neighborhood of $\partial W_0$ in $W_0$
and identify it with $\partial W_0 \times [0, \epsilon] \equiv
W_0(\epsilon) $, $\partial W_0  = \partial W_0 \times \{0 \}$. Then  for
any $\tau \in [0, \epsilon]$, we have the following: $$ W_0 \cap
\mathbb{R}^{N+1} \times [0, \epsilon] = W_0(\epsilon), \quad W_0(\epsilon)
\cap \partial W_0 \times \tau \subset \mathbb{R}^{N+1} \times \tau.$$
(2)\qua For a smaller $\epsilon' < \epsilon$, $p(\{ W_0(\epsilon') \cap
\partial W_0 \times \tau \})
\subset \mathbb{R}^{N+1}$ is independent of $\tau \in [0, \epsilon']$, where 
$p\co \mathbb{R}^{N+2} \to \mathbb{R}^{N+1}$ is the obvious projection.

(3)\qua $I|W_0 \backslash W_0(\epsilon)$ is a $C^{\infty}$ embedding
(ie without corners).

(4)\qua  Let $\overline{M}_1$ and $\overline{N}_0$ be tubular neighborhoods
of $M_1$ and $N_0$ in $\partial W_0$ respectively. As before we identify
small neighborhoods of $\partial \overline{M}_1$ and $\partial
\overline{N}_0$ with $\partial M_1 \times (0, \delta] \equiv M_1(\delta)$
and $\partial N_0 \times (0,\delta] \equiv N_0(\delta)$ respectively.
Then we have the following:

\quad (a)\qua $\overline{M}_1$ and $\overline{N}_0$ are diffeomorphic to $M_1$ and $N_0$ respectively,

\quad (b)\qua $\partial W_0 \cap \mathbb{R}^N \times \{0\} = M_1, \quad
\partial W_0 \cap \mathbb{R}^N \times \{1\} = N_0$,

\quad (c)\qua for $0< \tau < \delta$,  $M_1(\delta) \cap \mathbb{R}^N \times \tau$
and $N_0(\delta) \cap \mathbb{R}^N \times \{1- \tau \}$ are diffeomorphic to
$\partial M_1$ and $\partial N_0$ respectively.

(5)\qua $W_0(\epsilon) \cap \mathbb{R}^N \times \tau \times [0,\epsilon]$
and $W_0(\epsilon) \cap \mathbb{R}^N \times \{ 1-\tau \} \times [0,
\epsilon]$, $\tau \in (0, \delta]$, are diffeomorphic to $\partial M_1
\times [0, \epsilon] \cup M_1$ and $\partial N_0 \times [0,\epsilon]
\cup N_0$ respectively.

Moreover $W_0(\epsilon) \cap \mathbb{R}^N \times \tau \times [0,
\epsilon] \subset \mathbb{R}^N \times \{0\} \times [0, \epsilon]$ and
$W_0(\epsilon) \cap \mathbb{R}^N \times \{1- \tau \} \times [0, \epsilon]
\subset \mathbb{R}^N \times \{0\} \times [0, \epsilon]$ are independent
of $\tau $.

(6)\qua Let us denote $\overline{M}_1(\tau) = W_0(\epsilon) \cap \mathbb{R}^N \times \{\tau \}  \times [0, \epsilon]
= \{ (a,b,c) \in \overline{M}_1(\delta) \cap \mathbb{R}^N \times \mathbb{R} \times \mathbb{R}\}$.
We have similar notations for $N_0$. 
Then we have $\overline{N}_0(\tau) = \overline{M}_1^{\sharp}(\tau)  \equiv
\{(a, 1-b,c) | (a,b,c) \in \overline{M}_1(\tau)\}$, $\tau \in (0, \delta]$.

Let us denote $\widetilde{W}'_0= W_0 \backslash W_0 \cap \mathbb{R}^N \times 
\{ 0,1 \} \times [0, \epsilon]$.
Then $I|\widetilde{W}'_0$ is a smooth embedding with corners.
Let us denote by $g_0$ the Riemannian metric on $\widetilde{W}'_0$
induced from the standard one on $\mathbb{R}^{N+2}$ through $I$.
Notice important properties of $g_0$; for $\epsilon' < \epsilon$:
$$g_0|W_0(\epsilon') = (g_0|\partial W_0) \times dy^2, \quad
g_0|\overline{M}_1((0, \delta]) \cap W_0(\epsilon')  
 = (g_0|\partial M_1) \times dx^2 \times dy^2$$
where we have the coordinates $(t,x,y) \in \mathbb{R}^N \times \mathbb{R} \times
\mathbb{R}$. We have a similar expression for $g_0|\overline{N}_0(\delta)$.

Now we define the desired Riemannian metric on $\interior W_0$ as follows.
Let us put:
\begin{multline*}
\widetilde{W}_0 = \widetilde{W}_0' ~\cup~ [ \partial W_0 \backslash
  \{M_1 \cup N_0 \}] \times (-\infty, 0]  ~\cup~ \\
[ \partial M_1 \times ( - \infty, \epsilon] ~\cup~ M_1 \times \{ \epsilon \}
  ] \times (- \infty,0] ~\cup~ \\
[ \partial N_0 \times ( - \infty, \epsilon] ~\cup~ N_0 \times \{ \epsilon \}
  ] \times [1, \infty).
\end{multline*}
Clearly $\widetilde{W}_0$ is diffeomorphic to $\interior W_0$. On the
other hand, by the above remark, there is a natural extension of $g_0$
to $g|\widetilde{W}_0$. This is the desired complete Riemannian metric.
Later we will use the following notations. Let us choose large $t \gg 0$
and $s \gg 0$.  Then we denote:
\begin{multline}
\widetilde{W}_0(t) = [\partial W_0\backslash\{M_1 \cup N_0\}] \times \{-t\} \\
~\cup~ \partial M_1 \times \{ -t \} \times (- \infty, 0] \\
~\cup~ \partial N_0  \times \{ -t \} \times [1, \infty),
\end{multline}
$$M_1(s) =  \{ \partial M_1 \times ( - \infty, \epsilon] \cup M_1 \times \{ \epsilon \} \} \times \{-s\}.$$
We denote $M_1([S_0,S_1]) = \cup_{s \in [S_0,S_1]} M_1(s)$.
We have a similar notation $N_0(s)$. 

\subsection{Casson handles as Riemannian manifolds}
Let us generalize the previous construction over the triple
$(W_0,M_1,N_0)$ as follows.
Suppose we have $n+1$ submanifolds $M_0, M_1, \dots, M_n$ in $\partial W_0$
such that they are disjoint from each other, $M_i \cap M_j =\emptyset$.
Then by choosing an appropriate embedding
$I\co W_0 \hookrightarrow \mathbb{R}^N$,
one can construct a complete Riemannian metric on int$W_0$ as 2.A.
Then one expresses this as $\widetilde{W}_n$:
$$\widetilde{W}_n = \widetilde{W}_0' \cup [ \partial W_0 \backslash
  \cup_i M_i ] \times (-\infty, 0]  \cup_i
  [ \partial M_i \times ( - \infty, \epsilon] \cup
  M_i \times \{ \epsilon \} ] \times (- \infty,0].$$
Let $\Sigma$ be a Riemann surface of genus $g \geq 1$, and
$\overline{\Sigma}$ be its solid Riemann surface.
Thus $\overline{\Sigma}$ is a three--manifold with boundary $\Sigma$.
Let us put $W_g = \overline{\Sigma} \times [0,1]$, and take $g$ embeddings
of $S^1 \times D^2$ into $\partial W_g$
which represents a generating set of $\pi_1(W_g)$.
We denote their images by $M_1, \dots, M_g$.
One may assume that all $M_i$ and $M_j$ are disjoint.
A kinky handle with $k$ kinks is expressed as
$(W_k, N_0, M_1, \dots, M_k)$, where $N_0$ is another embedding
of $S^1 \times D^2$ into $\partial W_k$.
The isotopy type of $N_0$ is determined by the sign
$(\pm, \dots \pm) $ ($k$ times).  There are exactly $k+1$ diffeomorphism 
types of kinky handles with $k$ kinks, since the diffeomorphism type of the
kinky handle is invariant under permutation of the signs.

A Casson handle is constructed by taking end connected sums of kinky
handles infinitely many times. It is expressed by an infinite tree
(embedded in $\mathbb{R}^2$) with additional information.
 Let us introduce the following:
\begin{defn}
Let $T$ be an infinite tree with only  one end point $*$.
$T$ is said to be a signed tree if every edge $e$ is
assigned  one of $+$ or $-$. 
\end{defn}
Let $n(v)+1$ be an integer which is the number of edges of $T$ 
with common end-point $v$. A signed tree $T$ is said to be of
\definition{bounded multiplicity} if there exists a constant $C$ with
$C \geq n(v)$ for any $v \in T$.

Any tree with one end-point $*$ admits a natural distance from $*$.
One will denote it by $| \quad |$.
Let $T$ be a signed tree  with one end-point $*$. 
Let us take any vertex $v \in T$. Then there exists a unique 
edge $e $ on $T$ such that one of the end-point of $e$ is $v$, and the other
$v'$ satisfies $|v'| < |v|$. To specify this $e$, one denotes it as $e(v)$.
Let us denote finite subtrees by $T_v = \cup_{j=0}^{n(v)} e_j$, 
where all  the edges  $e_j$ has $v$ as one end-point, and
$e_0 =e(v)$. 
Recall that each $e_j$ is assigned with $\pm$.
Since   $e_j$, $1 \leq j \leq n(v)$,  are all assigned  one of $\pm$, 
the set $ \{ e_1, \dots, e_{n(v)} \} $  
determines an isotopy class of smooth embedding 
$N_0 \cong  S^1 \times D^2  \hookrightarrow \partial W_{n(v)}$,
where $W_{n(v)}$ is diffeomorphic to $\overline{\Sigma}_{n(v)} \times
[0, 1]$ and $e(v)$ corresponds to $N_0$.
This class determines a diffeomorphism type of a kinky handle.
Thus every $T_v$ corresponds to data $(W_{n(v)}, N_0,M_1, \dots, M_{n(v)})$,
a kinky handle with $n(v)$ kinks as described above.
In particular one may assign a complete Riemannian manifold $\widetilde{W}_v$
for every vertex $v \in T$.

Let us take an edge $e$ with end-points $v$ and $v'$,
and assign  a large number $S \gg 0$.
Then one can make end-connected sums of $\widetilde{W}_v$ with
$\widetilde{W}_{v'}$ as:
$$\widetilde{W}_{v'} \natural_{e, S} \widetilde{W}_v
\equiv \widetilde{W}_{v'} \backslash M_e([S, \infty)) \natural
\widetilde{W}_v \backslash N_v([S, \infty)).$$
One will also denote 
the set $(T,  S )$ by $T$, and also call it  a \definition{signed tree}.
For every $e$, one may connect $\widetilde{W}_v$ and $\widetilde{W}_{v'}$ as
above. By iterating these end-connected sums for every edge, 
one gets a complete Riemannian four--manifold $CH(T)$.  The diffeomorphism
type of this space, is the Casson handle corresponding to $T$.

We denote $(W_0, N_0,M_1)$ as the simplest kinky handle.
Then $\widetilde{Y}(S)$ is a periodic  Riemannian four--space:
$$ \widetilde{Y}(S)  = \cup_{j \in \mathbb{Z}} \widetilde{W}^j_0([-S,S]), \quad
  CH(\mathbb{R}_+) = \widetilde{Y}(S)_0 =
  \cup_{j \in \mathbb{N}} \widetilde{W}^j_0([-S,S])$$
where $\widetilde{W}_0^j([-S,S]) = \widetilde{W}_0^j \backslash
N_0([S, \infty)) \cup M_1([S, \infty))$,
 and 
one enumerates the same $\widetilde{W}_0$ by $\widetilde{W}_0^j$,
$j \in \mathbb{Z}$.  Even though $CH(\mathbb{R}_+)$ has very complicated
smooth structure, $\widetilde{Y}(S)$ is diffeomorphic to $\mathbb{R}^4$,
and so it would not  appropriate to write it as $CH(\mathbb{R})$.

Let $(W_l, N_0, M_1, \dots, M_l)$ be a kinky handle with $l$ kinks.
Then one defines the following another Riemannian manifolds.
Let us choose any  $1 \leq j \leq l$, say $j =1$.
Now one defines:
\begin{equation}
\begin{aligned}
 & V(S,l)_0 = \widetilde{W}_l \backslash N_0([S,  \infty)) \cup_{k=1}^l M_k([S, \infty)), \\
 & Y(S,l,1)_0 = V(S,l)_0 \cup_{j=2}^l \cup_{M_j(S) \sim N_0(S)} 
\widetilde{Y}(S)_0 \backslash N_0([S,  \infty)), \\
& Y(S,l,1) = Y(S,l,1)_0  / N_0(S) \sim M_1(S), \\
& \widetilde{Y}(S,l,1)_0 = \cup_{j \in \mathbb{N}}
\cup_{M_1^j(S) \sim N_0^{j+1}(S)} Y(S,l,1)_0^j, \\
& \widetilde{Y}(S,l,1) = \cup_{j \in \mathbb{Z}}
\cup_{M_1^j(S) \sim N_0^{j+1}(S)} Y(S,l,1)_0^j.
\end{aligned}
\end{equation}
There is an infinite signed tree with a base point $(T_{l,1})_0$
such that   $\widetilde{Y}(S,l,1)_0$ is diffeomorphic to
 $CH((T_{l,1})_0)$.

Let us take $l$ half-lines $\mathbb{R}_+, \dots \mathbb{R}_+$ ($l$ times) where
on each $\mathbb{R}_+$ the same sign is given.
Signs on two different half-lines may be mutually different.
Then one identifies all $\mathbb{R}_+$ at $0$. 
The result is a connected signed infinite
tree with a  base point. One denotes it by  $T_0^l$. 
$(T_{l,1})_0$ can be expressed as
$\mathbb{R}_+ \cup_{j \in \mathbb{N}} (T_0^{l-1})^j$,
where $(T_0^{l-1})^j$ are the same $ (T_0^{l-1})$ indexed by $j \in
\mathbb{N}$,
and  on $\mathbb{R}_+$ the same sign is given.

Let us define:
\begin{equation}
\begin{aligned}
& Y(S,l,m,1)_0 \equiv
        V(S,l)_0  \cup_{j=2}^l \cup_{M_1^j(S) \sim N^1_0(S)}  \widetilde{Y}(S,m,1)_0, \\
& Y(S,l,m,1) = Y(S,l,m,1)_0  / N_0(S) \sim M_1(S), \\
& \widetilde{Y}(S,l,m,1)_0 = \cup_{j \in \mathbb{N}}
\cup_{M_1^j(S) \sim N_0^{j+1}(S)} Y(S,l,m,1)_0^j, \\
& \widetilde{Y}(S,l,m,1) = \cup_{j \in \mathbb{Z}}
\cup_{M_1^j(S) \sim N_0^{j+1}(S)} Y(S,l,m,1)_0^j.
\end{aligned}
\end{equation}
Similarly one has $\widetilde{Y}(S,l,m,1)_0 = CH((T_{l,m,1})_0)$.
Let us take  signed infinite trees 
$(T_{m,1})_0, \dots, (T_{m,1})_0$ ($l$ times).
Then one identifies all $ (T_{m,1})_0$  at $0$. The result is a connected
signed infinite tree with a  base point. One denotes it by  $ (T_{m,1})_0^l$. 
As above, $(T_{l,m, 1})_0$ can be expressed as 
$\mathbb{R}_+ \cup_{j \in \mathbb{N}}  [(T_{m,1})_0^{l-1}]^j $,
where $   [(T_{m,1})_0^{l-1}]^j $  are the same    $  (T_{m,1})_0^{l-1}$  indexed by $j \in \mathbb{N}$,  and on $\mathbb{R}_+$ the same sign is given.
Similarly one uses a notation $(T_{l,m,1})_0^k$, which one identifies $k$  
$(T_{l,m,1})_0$  at $0$.
Then using this, one can obtain $(T_{k, l,m,1})_0$
by a similar method.

By iteration, one obtains 
$Y(S, l,m,n,1)$ and $\widetilde{Y}(S, l,m,n,1)_0 = CH((T_{l,m,n,1})_0)$
as above.
Inductively one can obtain $CH((T_{n_1, \dots, n_k,1})_0)$
by iteration.
One calls $ (T_{n_1, \dots, n_k,1})_0$  a \definition{homogeneous tree
of bounded type}.

Let $\overline{n} = \{ n_1, n_2, \dots, n_l, \dots\} $ be an infinite
set of positive integers.
One can iterate the previous construction infinitely many times. Then
one gets a complete Riemannian manifold:
$$CH((T_{\overline{n}})_0) = CH((T_{n_1, n_2, \dots})_0)$$
where we call $(T_{\overline{n}})_0$  a homogeneous tree.
If there is a bound $C$ with $n_j \leq C$ for all $j$,
then we call it a \definition{homogeneous tree of bounded multiplicity}.

\section{Some  properties of  elliptic operators over $\widetilde{W}_0$}
\subsection{Spectral decomposition}
Throughout this paper, one uses the following lemmas.
Let $P$ be an order--$1$ elliptic differential operator over a complete
Riemannian manifold
$W$ of bounded geometry. Let $w$ be a weight function over $W$ and
weighted Sobolev $k$--norms on $W$ by
$|u|W^k_w = (\Sigma_{l \leq k} \int \exp(w)|\nabla^l u|^2 vol )^{\frac{1}{2}}$.
Thus one has a bounded map $P\co W^{k+1}_w(W) \to W^k_w(W)$.
Let us denote by $P^*_w$ the formal adjoint operator with
respect to the weighted $L^2_w$ inner product.
\begin{lem} 
Suppose $\Spec\Delta = P^*_w \circ P $ is contained in $[\lambda^2, \infty)$ 
for some $\lambda >0$.
Then  $P\co W^{k+1}_w  \to W^k_w$ has closed range.

Similarly $P^*_w\co W^{k+1}_w \to W^k_w$ also has closed range.
\end{lem}
\begin{proof}
Let us consider the first statement.
One verifies this by induction. Let $k=0$ and
let $\{ u_i\}$ be a sequence in $W^1_w$ with $|P(u_i)|L^2_w \to 0$.
Since one has the estimate $|P(u)|L^2_w \geq \lambda |u|L^2_w$,
one sees $|u_i|L^2_w \to 0$.
Moreover by the elliptic estimate:
$$|u|W^{k+1}_w \leq C_k(|P(u)|W^k_w + |u|W^k_w), \quad k \geq 0$$
for some constants $C_k$, one concludes $|u_i|W^1_w \to 0$.
This shows the result for $k=0$.

We show by induction that if a sequence $\{u_i\}$  in $W^{k+1}_w$ 
satisfies $|P(u_i)|W^k_w \to 0$, then $|u_i|W^k_w \to 0$.
Then  by the elliptic estimate, one gets the conclusion.
Suppose the result holds for $k \leq k_0$, and take a sequence $\{u_i\}$  in $W^{k_0+2}_w$
with $|P(u_i)|W^{k_0+1}_w \to 0$.
Then since $|P(u_i)|W^{k_0}_w$ also converges to $0$, one knows that $|u_i|W^{k_0}_w$ 
converges to $0$
by the induction hypothesis.
By the elliptic estimate, one gets the result for $k_0+1$.
This completes the proof of the first statement.

Next let us consider the second statement. Let $H^{\perp} \subset L^2_w$
be a closed subspace which consists of the orthogonal complement of
$H \equiv P(W^1_w) \subset L^2_w$.  Notice the following relation:
$$\Delta(W^{k+2}_w) = P^*_w(H \cap W^{k+1}_w) \subset W^k_w.$$

Now let $\{ v_i\}$ be a sequence in $H \cap W^1_w$ with $|P^*_w(v_i)|L^2_w
\to 0$.  Since one may express $v_i = P(u_i)$ for $u_i \in W^2_w$,
one has $|\Delta(u_i)|L^2_w \to 0$. Thus one has $|u_i|L^2_w \to 0$.
Then by Cauchy--Schwartz, one has the estimate $|P(u_i)|^2L^2_w
\leq |\Delta(u_i)|L^2_w |u_i|L^2_w$.  This shows $|v_i|L^2_w \to
0$. By the elliptic estimate, one sees $|v_i|W^1_w \to 0$.  Then by
proceeding similarly as above, one verifies the conclusion by induction.
This completes the proof.
\end{proof}

Let $W$ and $w$ be as above, and consider an elliptic complex:
$$0  \longrightarrow   E_0   \stackrel{ d_0 }{\longrightarrow}
E_1  \stackrel{ d_1 }{\longrightarrow}  E_2  \longrightarrow  0$$
over vector bundles on $W$.
\begin{lem}
{\rm(1)}\qua Suppose $d_0\co W^{k+1}_w(E_0) \to W^k_w(E_1)$
has closed range, and $\ker d_0 =0$.
Then there exists a constant $C$ such that:
$$|(d_0)^*_w \circ d_0(u)| W^{k}_w(E_0) \geq C|u| W^{k+2}_w(E_0)$$
for any $u$. In particular $(d_0)^*_w \circ d_0$ gives an isomorphism.

{\rm(2)}\qua Suppose the corresponding bounded complex:
$$0  \longrightarrow  W^{k+2}_w(E_0)  \stackrel{ d_0 }{\longrightarrow}
  W^{k+1}_w(E_1)  \stackrel{ d_1 }{\longrightarrow}   W^k_w(E_2)
  \longrightarrow  0 $$
satisfies the following; 
(1) $d_0$ and $d_1$ have closed range, (2) the first cohomology group
$\ker d_1 / \im d_0 = 0$, and (3) $\ker d_0 = 0$.

Then for $P = (d_0)^*_w \oplus d_1$ and $\Delta = P^*_w \circ P$,
one has $\Spec\Delta \subset [\lambda, \infty)$ for some positive constant
$\lambda >0$.  Moreover $(d_1)^*_w$ has closed range.
\end{lem}
\begin{proof}
(1)\qua  Suppose the contrary, and take a sequence $\{u_i \} \subset W^2_w$
with  $|(d_0)^*_w \circ d_0(u_i)|L^2_w \to 0$. Then by Cauchy--Schwartz
one has a convergent sequence $|d_0(u_i)|L^2_w \to 0$.
Then by the elliptic estimate, one has $|d_0(u_i)|W^1_w \to 0$.
By closedness and the open mapping theorem, this shows existence of
non-trivial kernel.  This is a contradiction. The rest of the proof
follows by induction.

Part (2) follows from (1) and lemma 3.1. This completes the proof.
\end{proof}

Let us use the above notations. Let us put
$H = L^2_w(E_1) \cap (\ker d_1)^{\perp}$. 
Let us consider restriction of $\Delta =  (d_1)^*_w  \circ  d_1 |H$.
\begin{cor}
Suppose $d_0$ and  $d_1$ has closed range as above.  Then
$\Spec\Delta|H$ is contained in $[\lambda, \infty)$ for some positive
$\lambda >0$.
\end{cor}

\subsection{Analysis over $\widetilde{W}_0$}
Let us take the complete Riemannian manifold  $\widetilde{W}_0$ in 2.A.
Let us take a smooth function $w\co \widetilde{W}_0 \to [0, \infty)$.
Then as before one has weighted Sobolev $k$ norms on $\widetilde{W}_0$ by
$|u|W^k_w = (\Sigma_{l \leq k} \int \exp(w)|\nabla^l u|^2 vol )^{\frac{1}{2}}$.
Then we have the Atiyah--Hitchin--Singer complex (AHS complex) as:
$$0 \longrightarrow  W^{k+2}_w(\widetilde{W}_0; \Lambda^0)
  \stackrel{ d }{\longrightarrow}  W^{k+1}_w(\widetilde{W}_0; \Lambda^1)
  \stackrel{ d_+ }{\longrightarrow} W^k_w(\widetilde{W}_0; \Lambda^2_+)
  \longrightarrow  0$$
where $d^+= (1+*) \circ d$.

Let $M$ be a  complete Riemannian $3$ manifold without boundary,
but possibly  be non-compact.
We have in mind $M = \widetilde{W}_0(t), N_0(s)$ or $M_1(s)$.
Let $g$ be a Riemannian metric on $M$, and denote also by $g$ the 
product metric on $M \times \mathbb{R}$. 
Let $\Lambda^1(M \times \mathbb{R})$ and $\Lambda^2_+(M \times \mathbb{R})$ be
the exterior differentials on $1$ forms and self-dual $2$--forms with
respect to $g$.  Then we have the natural identification:
$$\Lambda^1(M \times \mathbb{R}) = p^*(\Lambda^1(M)) \oplus p^*(\Lambda^0(M)),$$
$$\Lambda^2_+(M \times \mathbb{R}) = p^*(\Lambda^1(M))$$
where $p\co M \times \mathbb{R} \to M$ is the projection. The isomorphisms are
given by:
$$u+v\,dt \leftrightarrow (u,v), \quad *_M u+ u \wedge dt \leftrightarrow u.$$
Let us put $X= M \times \mathbb{R}$.
Using the $L^2$ adjoint operator (here we do not use any weight
functions), we get the  elliptic operator
$P= d^*  \oplus  d^+\co  \Lambda^1(X) \to \Lambda^0(X) \oplus \Lambda^2_+(X)$.
Using  the above identification, this is expressed as:
$$P= d^*  \oplus  d^+\co p^*(\Lambda^0(M) \oplus \Lambda^1(M))
\to  p^*(\Lambda^0(M) \oplus \Lambda^1(M)).$$
Let us use $t$ as the coordinate on $\mathbb{R}$. 
Then by a straightforward calculation, one has the following expression:
$$P= - \frac{d}{dt} + \begin{pmatrix}
  *_M d &d \\
d^* & 0
\end{pmatrix}
\equiv - \frac{d}{dt} + Q.$$
Notice that $Q$ is an elliptic self adjoint  differential operator
on $L^2(M ; \Lambda^0(M) \oplus \Lambda^1(M))$.
Notice that if $M$ is closed, then $Q$ is not invertible (one has
constant functions). When $M$ is non compact, $0$ may be contained 
in the spectrum of $Q$, even if there are no kernels for $Q$.

Let us  introduce a weight function $\tau$ on $X$ as follows.
Let $M$ be possibly non-compact as above. Let us choose a positive number
$\delta >0 $ which will be specified later,  and choose a smooth  proper
function $\tau'\co \mathbb{R} \to [0, \infty )$ with
$\tau'| \mathbb{R} \backslash [-1,1](t) = \delta|t|$.
Then we define $\tau\co M \times \mathbb{R} \to [0, \infty)$ by $\tau(m,t) \equiv \tau'(t)$. Let us introduce
weighted Sobolev $k$ norms on $X$ by:
$$|u|W^k_{\tau} = \left(\Sigma_{l \leq k}
  \int \exp(\tau)|\nabla^l u|^2 \vol \right)^{\frac{1}{2}}.$$
 We denote these Sobolev weighted spaces by $W^k_{\tau}$ or $L^2_{\tau}$ (for $k=0$).
Then one  defines an isometry:
$$I\co L^2(X, p^*(\Lambda^0 \oplus \Lambda^1))
\to L^2_{\tau}(X, p^*(\Lambda^0 \oplus \Lambda^1))$$
by $I(u)= \exp(-\frac{\tau}{2})u$. Then one has the equality:
$$I^{-1}PI = P+ \frac{1}{2}\frac{d \tau}{dt}\co 
W^{k+1}(X;  p^*(\Lambda^0 \oplus \Lambda^1)) \to 
W^k(X;   p^*(\Lambda^0 \oplus \Lambda^1))    .$$
Let $d^*_{\tau}$ be the $L^2_{\tau}$ adjoint operator, i.e. 
$\langle u, d(v)\rangle |L^2_{\tau} =
  \langle d^*_{\tau}(u), v\rangle |L^2_{\tau}$,
and put
$P_{\tau} = d^*_{\tau} \oplus d^+$. Then we have the following expression:
$$I^{-1} P_{\tau} I = - \frac{d}{dt} + \begin{pmatrix}
  *_M d & d \\
d^* & -\frac{d \tau}{dt}
\end{pmatrix}
+ \frac{1}{2}\frac{d \tau}{dt} \equiv   
 - \frac{d}{dt} + Q_{\tau}.$$

We recall the definition of $\widetilde{W}_0(t)$ or $M_1(s),N_0(s)$
in the last paragraph in 2.A. Practically we will apply the next lemma
for  $M = \widetilde{W}_0(t)$.  Let us  consider the  operator:
$$Q_{\delta} = \begin{pmatrix}
*_M d  + \frac{\delta}{2}  & d \\
d^* & -\frac{\delta}{2}
\end{pmatrix}\co W^{k+1}(M  , \Lambda^1  \oplus  \Lambda^0 )
\to W^k(M  , \Lambda^1 \oplus  \Lambda^0  ).$$
Notice that $Q_{\delta}$ coincides with $Q_{\tau}$ above for $ t \geq 1$.
\begin{lem} Suppose $M$ is a closed Riemannian manifold.
Then for a small $\delta >0$ and $k \geq 0$,
$Q_{\delta}$
gives an isomorphism.  
\end{lem}
\begin{proof}
Let us take $(u,v) \in W^{k+1}(M, \Lambda^1 \oplus   \Lambda^0 ) $, and put
$Q_{\delta}(u,v) = (x,y)$.  
 Then we have:
$$ \left(*d+ \frac{\delta}{2}\right)u +dv = x, \quad
 d^*u - \frac{\delta}{2}v =y.$$
 For convenience, we recall 
$*^2u = (-1)^{p(n-p)}u$, $d^*u = (-1)^{np+n+1}*d*u$ for a $p$--form $u$
on an $n$--dimensional manifold.  Then from the above, we have:
$$ d^**du + \frac{\delta}{2} d^* u + d^*dv = d^* x  
= \left(\left(\frac{\delta}{2}\right)^2 + d^*d\right)v + \frac{\delta}{2} y.$$
Thus we can solve $v$ from $(x,y)$ as follows:
$$v= \left(\left(\frac{\delta}{2}\right)^2 + d^*d\right)^{-1}
  \left(d^*x + \frac{\delta}{2}y\right).$$
Next we have the following equality:
$$*^2du + \frac{\delta}{2}*u + *dv = *x 
= du + \frac{\delta}{2}*u + *dv.$$
\begin{equation}
\begin{aligned}
d^*du + \frac{\delta}{2}d^**u + & d^**dv = d^**x \\
 & =  d^*du + \frac{\delta}{2}* du 
= d^*du + \frac{\delta}{2}\left(x- \frac{\delta}{2}u-dv\right).
\end{aligned}
\end{equation}
Then from this, we have the following equation:
\begin{equation}
\begin{aligned}
\left(d^*d - \left(\frac{\delta}{2}\right)^2\right)u =
  & -\frac{\delta}{2}(x -dv) + d^**x, \\
dd^*u= &  \frac{\delta}{2}dv + dy.
\end{aligned}
\end{equation}
Thus one can again solve $u$ as:
$$u= \left(\Delta - \left(\frac{\delta}{2}\right)^2\right)^{-1}
  \left(dy +d^**x- \frac{\delta}{2}x+ \delta\,dv\right).$$
It follows  that there is a constant $C= C(M)$ such that
$\Spec |Q_{\delta}| \geq C(\frac{\delta}{2})^4$. Notice the following
inequalities:
\begin{equation}
 \begin{aligned}
 |(u,v)|W^{k+1} \leq & C\left(\frac{\delta}{2}\right)^{-2}\left|
  \left(\Delta \pm \left(\frac{\delta}{2}\right)^2 \right)(u,v)\right|W^{k-1} \\
\leq &  C\left(\frac{\delta}{2}\right)^{-4}\left(|x|W^k+|y|W^k\right)
  \leq C\left(\frac{\delta}{2}\right)^{-4}|Q_{\delta}(u,v)|W^k.
 \end{aligned}
\end{equation}
In fact one may assume that $C$
only depends on 
$\sup_{0 \leq l \leq 2}\sup_{x \in M}|\nabla^l R|(x)$, 
where $R$ is the curvature form with respect to
the Levi--Civita connection.  This completes the proof.
\end{proof}

Let $(P_{\tau})^*_{\tau}$ be the formal adjoint of the Sobolev weighted spaces. Then one considers
the same analysis for   $(P_{\tau})^*_{\tau}$. Notice the equality: 
$$(I^{-1}  \circ P_{\tau} \circ I)^* =
  I^{-1} \circ  (P_{\tau})^*_{\tau}  \circ I.$$
Thus one can analyze an operator of the form
$P_{\tau}^*= \frac{d}{dt} + Q_{\delta}$, since $Q_{\delta}$ is
self-adjoint over $M$.

Let us take a smooth (non-proper) function:
$$\tau\co \widetilde{W}_0 \to [0, \infty)$$
which is horizontally constant, namely for some 
positive constant $\delta >0$ we have $\tau|\widetilde{W}_0(t) = \delta t$
for $t \geq 0$.
Using $\tau$, one defines weighted Sobolev spaces 
$W^k_{\tau}((\widetilde{W}_0, g), \Lambda_*)$.
One has the expression $P_{\tau}|M_1([0, \infty)) = -\frac{d}{ds} + Q_{\tau'}$
where:
$$Q_{\tau'} =   \begin{pmatrix}
  *d & d \\
  d_{\tau'}^* & 0
  \end{pmatrix}
  \co W^{k+1}_{\tau'}( M_1(s); \Lambda^1 \oplus  \Lambda^0) \to
  W^k_{\tau'}(M_1(s); \Lambda^1 \oplus \Lambda^0)$$
is an elliptic operator over $M_1(s)$, and  $\tau' = \tau|M_1(s)$ is
its restriction.  It follows that $Q_{\tau'}$ has closed range with
empty kernel.  We do not use this fact later.  In this case one can not
see whether  $P_{\tau}$ defined above is Fredholm or not, since
 one has no simple statement
as above for the weighted adjoint operator $(Q_{\tau'})^*_{\tau'}|
W^k_{\tau'}(M_1(s); \Lambda^1)$.

Let us take  another weight function $\mu\co \widetilde{W}_0 \to [0, \infty)$,
with $\mu|M_1(s),N_0(s) \equiv \delta' s$, where we choose a sufficiently small
$0< \delta' \ll \delta$.
Let us put:
$$w = \tau + \mu\co\widetilde{W}_0 \to [0, \infty).$$
We call $w$ a \definition{weight function} with \definition{weight}
$(\delta, \delta')$.  Let us use  $w$ as the weight function.   
Then one obtains a bounded operator:
$$P_w = d^*_w \oplus d^+\co W^{k+1}_w(\widetilde{W}_0; \Lambda^1) \to 
W^k_w(\widetilde{W}_0; \Lambda^0 \oplus  \Lambda^2_+).$$
Passing through the isometry $I\co L^2 \to L^2_w$, $I(u) = \exp(-\frac{w}{2})u$,
 $P_w$ has the following
expression  on $M_1([0, \infty))$ and $N_0([0, \infty))$:
$$  I^{-1} \circ P_w \circ I =  - \frac{\partial}{\partial s} + \begin{pmatrix}
*d & d \\
d^*_{\tau'}  & 0
\end{pmatrix} 
+ \frac{\delta'}{2} \begin{pmatrix}
1 & 0 \\
0 & -1
\end{pmatrix} 
 \equiv  - \frac{\partial}{\partial s} + Q^v_w.$$
On the other hand  one has another expression on $\widetilde{W}_0([0, \infty))$:
$$  I^{-1} \circ P_w \circ I   =   - \frac{\partial}{\partial t} + 
\begin{pmatrix}
*d  + \frac{\delta}{2}& d \\
d^*_{\mu} & -\frac{\delta}{2}
\end{pmatrix} 
  =  - \frac{\partial}{\partial t} + Q_w^h  .$$
\subsection{Decay estimate}
Let $M$ be a (possibly non-compact) $n$--dimensional Riemannian
manifold  of bounded geometry. In practice $M$ stands for
$\widetilde{W}_0(t)$ (or $N_0(s), M_1(s)$).
Let us take a smooth function $\mu\co M \to [0, \infty) $.
Let $E$ and $F$  be vector bundles over $M$, 
and take an elliptic operator:
$$Q\co W^{k+1}_{\mu}(M;E) \to W^k_{\mu}(M;F).$$ 
Let us choose $\tau \co M \times \mathbb{R} \to \mathbb{R}_+$ by $\tau (m,t) = \delta |t|$ for
$|t| \geq 1$, and  put $w = \tau + \mu$.
Then one considers 
the following operator:
$$P = - \frac{d}{dt} + Q\co W^{k+1}_w( M  \times \mathbb{R}, p^* E)
\to W^k_w(M \times \mathbb{R} ; p^* F).$$

Let $P^*_w$ be the formal adjoint operator with respect to 
the weighted inner product. Then $ P^*_w P$ admits
a spectral decomposition. For $u \in W^k_w(M \times
\mathbb{R}, p^*(E))$, we denote $\Spec u \in [0, \epsilon]$ if
$u$ lies in the image of the spectral projection
${\mathfrak{P}}([0, \epsilon])$. 

Suppose the above $Q$ is a Fredholm operator, and 
is self adjoint with respect to the weighted $L^2_{\mu}$ inner product.
\begin{prop}{\rm\cite{lockhart and mcowen,floer,taubes2}}\qua
For sufficiently small $\delta >0$, the above $P$ is Fredholm.
Moreover if $Q$ is invertible, then one may choose $\delta =0$.
In particular $P$ is also invertible.

Suppose $Q\co  W^{k+1}_{\mu}(M;E) \to W^k_{\mu}(M;F)$ has closed range
with $\ker Q=0$.
Let us put $\Delta = Q^*_{\mu} \circ Q$, and choose a positive
constant $\lambda_0 >0$ with $\Spec \Delta \in [\lambda_0^2, \infty)$.
Then  under the  above  condition, one has the following;
for any small $\epsilon \ll \lambda_0$, 
there exist positive constants $C,  \delta_0 >0$ such that for any $u \in 
W^{k+l}(M \times \mathbb{R},p^*E) $ with $\Spec u \in [0, \epsilon]$
for $P^*_w \circ P$, we have the following decay estimate:
$$|u(~~, t)|L^2_{\mu}(M) \leq C \exp(- \delta_0 |t|)
|u|L^2_{\mu}(M \times [0, \infty)), \quad t \geq 0.$$
\end{prop}
\begin{proof}
The first statement follows from the decay estimate.
Next let us take $u$  as above, and denote its slice on 
$M \times   s$ by $u_s$.
Let us put a smooth function  $f\co \mathbb{R} \to [0, \infty) $ by
$f(s) = \int_{M_s} \exp(\mu) |u_s|^2$. Then by differentiating,
one has the inequalities:
\begin{equation}
\begin{aligned}
  \quad f''(s) & = 2\left[ \int_{M_s} \exp(\mu) \langle u_s', u_s'\rangle + 
       \int_{M_s} \exp(\mu) \langle u_s'' ,u_s\rangle\right] \\
& \geq  2\left[ \int_{M_s} \exp(\mu) \langle u_s', u_s' \rangle + 
       \int_{M_s} \exp(\mu) \langle (\Delta-\epsilon) u_s ,u_s\rangle\right] \\
& \geq      2 \int_{M_s} \exp(\mu) \langle( \Delta - \epsilon)u_s ,u_s\rangle 
\geq 2(\lambda_0 - \epsilon)f(s).
\end{aligned}
\end{equation}
Thus one has differential inequalities 
$f'' \geq \lambda_0 f$. From this, one can get the 
desired estimate.
This completes the proof.
\end{proof}

Let $Y$ be a complete Riemannian four--manifold such that except for a compact
subset $K \subset Y$, $Y \backslash K$ is isometric to $Y(0) \times
[0, \infty)$ where $Y(0)$ is a closed manifold. Let us equip a weight
function $\tau\co Y \to [0, \infty)$ by $\tau | Y(0) \times t = \delta t$.
Let $I\co L^2(Y) \cong L^2_{\tau}(Y)$ be the isometry.
\begin{cor}{\rm\cite{lockhart and mcowen}}\qua
The corresponding:
$$P = I^{-1} P_{\tau} I\co W^{k+1}(Y; \Lambda^1) \to W^k(Y; \Lambda^0
  \oplus \Lambda^2_+)$$
gives a Fredholm map.
\end{cor}
\begin{proof}
This is well known, but for convenience we recall the proof. Let us take
a cut-off function
$\varphi|Y$ with $\varphi|Y(0) \times [1, \infty) \equiv 1$,
and $\varphi|K \equiv 0$. Let us choose any $u \in W^{k+1}(Y; \Lambda^1)$. Then
one may regard $\varphi u \in W^{k+1}(Y(0) \times \mathbb{R}; \Lambda^1)$.
Thus one has the estimate:
\begin{equation}
\begin{aligned}
|u|W^{k+1} & \leq |\varphi u|W^{k+1} + |(1- \varphi)u|W^{k+1} \\
& \leq C\{ |\varphi P(u)|W^k +
|[\varphi, P]u|W^k + |(1- \varphi)Pu|W^k \\
& \quad + |[(1-\varphi),P]u|W^k + |(1- \varphi)u|W^k\}.
\end{aligned}
\end{equation}
This shows that $P$ has closed range with finite-dimensional kernel.
One has a similar estimate for $P^*$.
This shows the result.
\end{proof}

\subsection{First cohomology}
Let $w= \tau + \mu$ be as in the last paragraph of 3.B.
\begin{lem} For a sufficiently small choice of
$\delta' \ll \delta$,
$d\co W^{k+1}_w(\widetilde{W}_0; \Lambda^0)
\to W^k_w(\widetilde{W}_0; \Lambda^1)$ has closed range with $\ker d=0$.
\end{lem}
\begin{proof}
Clearly $\ker d=0$.

As lemma 3.3, passing through the isometry $I$,
$I^{-1} \circ d \circ I|\widetilde{W}_0([0, \infty))$ can be expressed as:
$$\frac{d}{dt} + \left(d, - \frac{\delta}{2}\right)\co
  W^{k+1}(\widetilde{W}_0([0, \infty)) \to
  W^k(\widetilde{W}_0([0, \infty); \Lambda^0 \oplus \Lambda^1).$$
Clearly this has closed range. One has the same property on $N_0([0, \infty))$
and $M_1([0, \infty))$.

Now let $H_i \subset W^{k+1}(\widetilde{W}_0)$, $i=0,1,2,3$ be closed
subsets whose supports lie in $N_0([0, \infty))$, $M_1([0, \infty))$,
$\widetilde{W}_0([0, \infty))$ and $K$ respectively,
where $K \subset \widetilde{W}_0$ is some compact subset. 
Then one may assume $\cup_i H_i = W^{k+1}(\widetilde{W}_0)$. Now by the 
open mapping theorem and above, every image $d(H_i)$
has closed range. This completes the proof.
\end{proof}

Let us consider $\widetilde{W}_0$ and the weight function $w$ on it 
with the weight constants $(\delta, \delta')$.
Let us recall that we have a bounded  complex:
$$0 \longrightarrow W^{k+2}_w(\widetilde{W}_0) \overset{d}{\longrightarrow} 
W^{k+1}_w(\widetilde{W}_0; \Lambda^1)
 \overset{d_+}{\longrightarrow} 
W^k_w(\widetilde{W}_0; \Lambda^2_+) \longrightarrow 0$$
where the weight function is $w= \tau + \mu$.

It is clear $\ker d_0 =0$.  Suppose 
$\omega \in W^1_w(\widetilde{W}_0; \Lambda^1)$ with $d_+( \omega ) =0$.
Then  by integration by parts, one has
$0= \int_{\widetilde{W}_0} |d_+\omega|^2 =
\frac{1}{2} \int_{\widetilde{W}_0} |d \omega|^2$.
In particular we have $d \omega=0$. The following sublemma shows $H^1(\AHS) =0$.

Let $W= \widetilde{W}_0$ or $\widetilde{Y}(S)$ in 2.B. Recall that
one has a weight function $\tau$ over $\widetilde{W}_0$. There is a
natural extension of $\tau$ on $\widetilde{Y}(S)$.
Let $\mu\co \widetilde{Y}(S) \to [0, \infty)$ be another weight function
defined as $\mu|\widetilde{W}^n_0([-S,S])(x) = \delta' (|n|+ t(x))$,
where $t\co \widetilde{W}^0_0([-S,S]) \to [0,1]$ satisfies
$t| N_0(S) \equiv 0$ and $t|M_1(S) \equiv 1$.

Let $w = \tau + \mu$ be a weight function on $W$. In order to make explicit the
weight constants, sometimes one uses notations $\tau(\delta)$ and
$\mu(\delta')$.  Later one will show that $d\co
W^{k+1}_w(\widetilde{Y}(S)) \to W^k_w(\widetilde{Y}(S); \Lambda^1)$ has
closed range. Assuming this, one has the following: 
\begin{sublem}
Let $f \in C^{\infty}(W)$ satisfies $|df|L^2_w(W) < \infty$.
Then there are constants $C$ and  $k= k(f)$ such that:
$$|f-k|L^2_w < \infty.$$
\end{sublem}
\begin{proof}
The situation differs from  \cite[lemma 5.2]{taubes},
in that on the end, we have non compact slices, $\widetilde{W}_0(t),
M_1(s)$, and $N_0(s)$.
In order to verify this, one takes two steps. First one sees a weaker
version of the sublemma.

Step 1: Let us choose any $\delta_0 < \delta$ and any $\delta_0' < \delta'$ 
so that the pair $(\delta_0, \delta_0')$ is
sufficiently near $(\delta , \delta')$.
Let us put $w= \tau(\delta)+ \mu(\delta')$ and
$w' = \tau(\delta_0) + \mu(\delta_0')$. Then we claim that $f- k \in L^2_{w'}$.

By lemma $3.2 (1)$ and lemma $3.4$, one has an estimate:
$$ \int_W \exp(w)|u|^2 \leq C \int_W \exp(w)|du|^2.$$
This shows that one has a Hilbert space, which is a completion of $C^{\infty}_{\text{cp}}(W)$
by the norm $|u|^2 = \int_W \exp(w)|du|^2$. 
Let us take any $f \in C^{\infty}(W)$ with $df \in L^2_w$.
Using a functional $s\co H \to \mathbb{R}$, by $s(v) = \int_W
\exp(w)(|dv|^2 - 2\langle dv, df\rangle)$,
one gets a unique $u = (d_w^* d)^{-1}(d_w^*df) \in H$ with $d_w^*d(f-u) =0$.
We show that $g = f-u$ is in $L^2_{w'}$.
 From $dg \in L^2_w$, one has estimates
$\int_{W([t,t+1])} |dg|^2 \leq \exp(- \delta t) \epsilon_t$,
where $\epsilon_t \to 0$. Notice that $dg$ satisfies the elliptic equation
$(d_w^* \oplus d)dg =0$. Thus  by the local Sobolev embedding and the above estimate, one gets:
$$\textstyle{\sup_{\widetilde{W}_0(t)}} |dg|
  \leq C \exp(- \frac{\delta t}{2})\epsilon_t.$$
In particular one has the next estimate:
$$|g(m,t)| \leq C_1 + \int_0^t C_2 \exp\left(-\frac{\delta s}{2}\right)
  \epsilon_s ds,$$
where $C_1 = |g(m,0)|$. This shows that there is a constant $C_m$ such that 
$\lim_{t \to \infty} \sup_{W(m,t)}|g- C_m|=0$. 
It is clear that $C_m$ is independent of $m$, and we write the same
constant by $C$.

Let us make a similar procedure on $s$ direction along  $M_1(s)$ and
$N_0(s)$. Then one finds other constants $C'$ and $C''$ such that
$g-C'|M_1(m,s)$ and $g-C''|N_0(m,s)$ vanish at infinity for every
$m$. Again it is clear that $C=C'=C''$.

Now one finds a constant $c$ with $|dg|(x) \leq c \exp(- \frac{w(x)}{2})$. 
Then by integration, one has the inequality:
$$|g- C|(x) \leq c \exp\left(- \frac{w(x)}{2}\right).$$
Notice that for any $\alpha >0$, $w^{\alpha}$ is integrable over $W$.
Combining with this and the above inequality, one gets the first claim.

Step $2$: For any sufficiently small pair $(\delta, \delta')$, one has
a well-defined first cohomology group $H^1(\delta, \delta') = \ker d
\cap L^2_w(W; \Lambda^1)/ d(W^1_w(W))$.
Then one has a natural map:
$$i\co H^1(\delta, \delta') \to H^1(\delta_0, \delta_0').$$
We claim that $i$ is an injection. This is enough for
the proof of the sublemma.

Let us put $\Delta(\delta, \delta')= (d_w^* \oplus d)_w^*(d_w^* \oplus d)$
on $L^2_w(W; \Lambda^1)$, $w= w(\delta, \delta')$. Then the spectrum of
$\Delta(\delta, \delta')$ is discrete near $0$.
Thus for $(\delta_0, \delta_0')$ sufficiently near $(\delta, \delta')$,
one has $\dim H^1(\delta, \delta') \geq \dim H^1(\delta_0, \delta_0')$.

Suppose these dimensions are different. Then one chooses another $\delta_1
< \delta_0$ and $\delta_1' < \delta_0'$. If $\dim H^1(\delta_0, \delta_0')
= \dim H^1(\delta_1, \delta_1')$, then one puts $\delta_0 = \delta$
and $\delta' = \delta_0'$, and gets the conclusion. If these dimensions
are different, one has $\dim H^1(\delta_0, \delta_0') >$ $\dim
H^1(\delta_1, \delta_1')$.
If they  are different, then one must have $\dim H^1(\delta, \delta') \geq 2$.
Let us take another $\delta_2 < \delta_1 < \delta_0$ and $\delta_2'$ similarly.

One iterates this process until one gets the equality on dimension.
If this is not the case, one has an infinite sequence $(\delta_i, \delta_i')$.
By choice, one may assume that this sequence converges to
$(\delta_{\infty}, \delta_{\infty}')$ with $2\delta_0(') - \delta(') <
\delta_{\infty}(') < \delta(')$.
Then $\Delta(\delta_{\infty}, \delta_{\infty}')$ must have a continuous
spectrum near $0$.  This is a contradiction.  This completes the proof.
\end{proof}

\begin{rem}
{\rm When $\delta'=0$ (no weight in the horizontal direction),
one can still find some constant $C$ such that $g -C$ vanishes at
infinity.  Let us take any $(m', t) \in \widetilde{W}^n_0([-S,S])
\subset \widetilde{Y}(S)$. Then we show that if $n$ is sufficiently large,
then $|g(m',t)-C_m|$ is sufficiently small for $C_m$ as above.
With respect to $n$, let us choose a sufficiently large $T$.  Then since
one has the estimate
$|dg|\widetilde{Y}(S)(T) \leq C \exp(- \frac{\delta T}{2})$, one may
assume $|g(m,T) - g(m', T)|$ is sufficiently small. On the other hand,
$|g(m',t) - g(m', T)|$ is less than $\epsilon_n$, where $\epsilon_n \to 0$
as $n \to \infty$. Thus $|g(m',t)- C_m|$ is sufficiently small.}
\end{rem}
\begin{lem}
Let $W= \widetilde{Y}(S)$. $H^1(\delta, 0)$ is naturally
isomorphic to $H^1(\delta, \delta')$.
\end{lem}
\begin{proof}
Let $\alpha \in H^1(\delta, 0)$ be a harmonic representative.
Then from corollary 4.1 below, and by the same method as \cite[lemma
5.3]{taubes},
one can verify $\alpha \in L^2_{w(\delta, \delta')}(\widetilde{Y}(S);
\Lambda^1)$.
This completes the proof.
\end{proof}

\section{Fourier--Laplace transforms between open manifolds}
\subsection{Fourier--Laplace transforms}
Our main application is  the analysis over $(W_0,N_0,M_1) =$ kinky
handles.  Recall that for the simplest kinky handle, $W_0$ is
diffeomorphic to $S^1 \times D^3$.  In this section, one considers $(W_0 =
S^1 \times D^3, M_1, N_0)$, where $M_1,N_0$ are all diffeomorphic to
$S^1 \times D^2$, and the embeddings are given as follows.
Let us regard $\partial W_0 = S^1 \times D^2_+ \cup S^1 \times D^2_-$,
and consider two knots $\mu \subset  S^1 \times D^2_+$, $\gamma \in S^1
\times D^2_-$ as follows; 

(1)\qua $\mu$ represents a Whitehead link diagram (see \cite[page 79, figure
3.3]{kirby}) and 

(2)\qua $\gamma = S^1 \times 0 \subset S^1 \times D^2_-$.

Notice that $\mu$ is null-homotopic, and $\gamma$ represents
a generator of $\pi_1(W_0)$.
Let us take tubular neighborhoods of $\mu$ and $\gamma$ in $S^1 \times
D^2_{\pm}$ respectively. We denote them as $\overline{\mu}$ and
$\overline{\gamma}$.
Then we choose $\overline{\mu}=M_1$ and $\overline{\gamma}=N_0$.
Let us denote the quotient space by $Y = W_0 / M_1  \sim N_0$.
By a simple calculation, one knows the following:
\begin{sublem} $H^1(Y; \mathbb{R}) =\mathbb{R}$ and
$H^2(Y;\mathbb{R}) =H^3(Y;\mathbb{R})=0$.
\end{sublem}
In this section, one introduces the Fourier--Laplace
transform over $Y$, which is an open analogue of the one in \cite{taubes}.
In order to do this, one will only require the above topological conditions.
In particular one can apply the  method in this section for  cylindrical
manifolds $Y$ which includes the simplest Casson handles.

Now by the previous construction, one has a complete Riemannian manifold
$\widetilde{W}_0$.
Let us take  infinitely many $\widetilde{W}_0$, and for convenience,
we enumerate the same $\widetilde{W}_0$ as $\widetilde{W}_0^i$, $i \in
\mathbb{Z}$.
Let us choose a sufficiently large $S$, and introduce  notations:
\begin{equation}
\begin{aligned}
 & \widetilde{W}_0([-S,S]) =
\widetilde{W}_0\backslash \{ M_1([S, \infty)) \cup N_0([S, \infty)) \}, \\
& Y(S) = \widetilde{W}_0([-S,S])/ N_0(S) \sim M_1(S).
\end{aligned}
\end{equation}
One has a natural extension of notation $Y(S)(t) = \widetilde{W}_0(t) /
N_0(S,t) \sim M_1(S,t)$, $t \in [0, \infty)$.
Recall that one has a weight function $\tau\co \widetilde{W}_0 \to [0, \infty)$
by $\tau|\widetilde{W}_0(t) = \delta t$.
One can make a natural extension of $\tau$ on $Y(S)$.
We will use the same $\tau$  for this.

Now one has a periodic Riemannian manifold as:
\begin{equation}
\begin{aligned}
\widetilde{Y}(S)=  & \dots \cup_{M_1^{-n}(S) \sim N_0^{-n+1}(S) } 
  \widetilde{W}_0^{-n+1}([-S,S]) \\
& \cup _{M_1^{-n+1}(S) \sim N_0^{-n+2}(S)} \widetilde{W}_0^{-n+2}([-S,S]) \dots
\end{aligned}
\end{equation}
There is a free and isometric $\mathbb{Z}$--action on $\widetilde{Y}(S)$, and  
we denote the action $1 \in \mathbb{Z}$  by $T$.

Let us take  vector bundles $E_i \to Y(S)$ (for $i=1,2$) and a differential
operator $D$ between them.
Then it lifts to the following one:
$$D\co C^{\infty}_{\text{cp}}(\widetilde{E}_1) \to
  C^{\infty}_{\text{cp}}(\widetilde{E}_2)$$
where $\widetilde{E}_i \to \widetilde{Y}(S)$ are the natural lifts.

Let us take any $\psi \in  C^{\infty}_{\text{cp}}(\widetilde{E}_1) $
and $z \in \mathbb{C}^*$.  Then we define the Fourier--Laplace transform
of $\psi$ by:
$$\widehat{\psi}_z(~~) = \Sigma_{n=-\infty}^{\infty} z^n(T ^n \psi)(~~).$$
Let  us define another vector bundle over $Y(S) \times \mathbb{C}^*$ as:
$$E_1'  \equiv [ \widetilde{E}_1 \otimes_{\mathbb{R}} \mathbb{C} ]/ \mathbb{Z} \to Y(S) \times \mathbb{C}^*$$
where $1 \in \mathbb{Z}$ sends $(\rho, \lambda) \in \widetilde{E}_1 \otimes_{\mathbb{R}} \mathbb{C}$
to $(T \rho, z\lambda)$. One may regard $E_1'$ as a family of vector bundles $\{ E_1'(z) \}$
over $Y(S)$, parameterized by $z \in \mathbb{C}^*$.

Now by restriction $\widehat{\psi}| \widetilde{W}_0^0( [-S,S])$,
$\widehat{\psi}_z$ defines a smooth section over $Y(S)$ of $E_1'(z)$,
where $\widetilde{W}_0^0 ([- S,S]) \subset \widetilde{Y}(S)$. Thus one
may regard $\widehat{\psi}_z \in C_{\text{cp}}^{\infty}(E_1'(z))$.

The Fourier--Laplace inversion formula is as follows.
Let us take  a smooth section $\widehat{\eta} \in
C^{\infty}_{\text{cp}}(E_1')$ with $\widehat{\eta}_z \in
C_{\text{cp}}^{\infty}(E_1'(z))$.  Then for any $s \in (0, \infty)$,
the following:
$$T^n \eta(x) \equiv \frac{1}{2 \pi i} \int_{|z| =s} z^{-n}
  \widehat{\eta}_z (\pi(x)) \frac{dz}{z}$$
defines a smooth section over $\widetilde{E}_1 \to \widetilde{Y}(S)$,
where
$\pi\co \widetilde{Y}(S) \to Y(S)$ is the projection
and
$x \in \widetilde{W}_0^0([-S,S])$.
By Cauchy's formula, these are inverses of each other, independently of $s$.

Let $D$ be as above. Passing through the Fourier--Laplace transform, one has
another differential operator between $(E_0)'$ and $(E_1)'$ by:
$$\widehat{D}_z \widehat{\psi}_z \equiv (\widehat{D \psi})_z.$$

Notice that at $z=1$, $(E_i)'$, $i=0,1$, are isomorphic to $E_i \to Y(S)$ respectively.
In fact every $(E_i)'(z)$ is isomorphic to $E_i$ as follows.
Let us take $t\co \widetilde{W}_0([-S,S]) \to [0,1]$, a smooth map such that
$t \equiv 0$ near $N_0(S)$, and $t \equiv 1$ near $M_1(S)$.
By taking  $[0, \infty) \subset \mathbb{C}$ as a branched line,
one may define $\log z$. Let us put $z^t \equiv \exp(t \log z)$, 
and consider $z^t \widehat{\psi}_z$. This gives a $\mathbb{C}$--valued
section over $E_0'$.  One calculates $\widehat{D}_z$ in terms of this
identification. The result is as:
$$\widehat{D}(z) = D + z^t[D, z^{-t}]\co C^{\infty}_{\text{cp}}(E_0) \to 
C^{\infty}_{\text{cp}}(E_1).$$
\subsection{Elliptic complexes over periodic covers}
Let $\{E_i, D_i \}_{i=0,1,2}$ be an elliptic complex over $Y(S)$.
It is clear that  it  gives a bounded one as:
$$0 \longrightarrow   W^{k+2}_{\tau}(\widetilde{E}_0 )
  \stackrel{ D_0 }{\longrightarrow}  W^{k+1}_{\tau}(\widetilde{E}_1)
  \stackrel{ D_1}{\longrightarrow}  W^k_{\tau}(\widetilde{E}_2)
  \longrightarrow  0.$$
\begin{prop}
The above $D_i$ is an acyclic Fredholm complex, ie all cohomologies
$H^i$ vanish, if  for all $z \in C =\{ w \in \mathbb{C}^*;
|w| = 1 \}$, $\{ \widehat{D}_z^i\co W^{k+i}_{\tau} 
 \to W^{k+i-1}_{\tau} \}$
are also so.
\end{prop}
\begin{proof}
For any $u \in W^k_{\tau}(\widetilde{Y}(S); \widetilde{E}_i)$,
$\widehat{u}_z$ has the following property; it lies in $W^k_{\tau}((E_i)_z
; Y(S)_z)$ for all $z$ with $0 < |z| \leq 1$.

When $|z|=1$, the inner product on $E_i$ has a natural extension
of a Hermitian metric on $E_i(z)$. 
Recall that $D_i\co W^{k+1}_{\tau}(Y(S); E_i(z)) \to W^k_{\tau}(Y(S);
E_{i+1}(z))$ has closed range.  Let $(D_i)^*_{\tau}$ be the formal
adjoint operator.
This differential operator is defined independently of $z$, with $|z|=1$
over $E_i(z)$.
By the assumption, one has an isomorphism:
$$D =  D_0 \oplus (D_1)^*_{\tau}\co
 W^{k+1}_{\tau}(Y(S); E_0(z) \oplus E_2(z)) \to 
W^k_{\tau}(Y(S);E_1(z)).$$
We denote by $D^{-1}$ its inverse.
On the other hand one has a natural lift of $D$ as:
$$D\co W^{k+1}_w(\widetilde{Y}(S) ; \widetilde{E}_0 \oplus
\widetilde{E}_2) \to W^k_w(\widetilde{Y}(S); \widetilde{E}_1).$$
Here we show that there is another bounded operator:
$$R\co W^k_w(\widetilde{Y}(S); \widetilde{E}_1) \to 
 W^{k+1}_w(\widetilde{Y}(S); \widetilde{E}_0 \oplus \widetilde{E}_2)$$
with $R \circ D $ is the identity.
 From this, one knows that $D$ has closed range. In particular
$D_0$ has closed range.

Next one proceeds similarly for 
$D^*_{\tau} = (D_0)^*_{\tau} \oplus D_1$.
Then one also knows that its lift $D^*_{\tau}$
has closed range. Clearly  it follows that $D_1$ has also closed range.

Now let us take $\psi \in W^{k+1}_w(\widetilde{Y}(S); \widetilde{E}_0
\oplus \widetilde{E}_2)$.
Then one considers $b(\psi) = D^{-1}( \widehat{\psi}) \in
W^{k+1}_{\tau}(Y(S); E_1(z))$.
This element also has the  property in the first paragraph of the proof.

Then one defines $R(\psi) \in W^{k+1}_w(\widetilde{Y}(S);\widetilde{E}_1)$ by:
$$T^n R(\psi)| \widetilde{W}_0^0([-S,S]) (x)= 
\frac{1}{2 \pi i} \int_C  z^{-n}( D^{-1}( \widehat{\psi}) (\pi(x))) \frac{dz}{z}.$$
We claim that this assignment:
$$R\co  W^k_w(\widetilde{Y}(S); \widetilde{E}_1) \to
  W^{k+1}_w(\widetilde{Y}(S); \widetilde{E}_0 \oplus \widetilde{E}_2)$$
gives a bounded operator.
Then by the definition of $D_i$ over $Y(S)$, it satisfies 
$R \circ D =1$.

Let us put $ C = \{ z =\exp(i \theta) \in \mathbb{C}  \}$, $\psi_n|\widetilde{W}^0_0([-S,S]) = T^n(\psi)$.
By regarding $b(\psi)_z \in W^{k+1}_{\tau}(\widetilde{W}^0_0([-S,S]))$,
one has a Fourier expansion
$b(\psi)= \Sigma_n z^n \varphi_n$, $\varphi_{n+1}|N_0(S) = \varphi_n|M_1(S)$.
Notice the equality $D(\varphi_n) = \psi_n$.
One can regard the value of the inner product 
$\langle b(\psi), b(\psi)\rangle |W^k_{\tau}(Y(S); E_*(z))$ as a
real-valued function of $\theta \in [0, 2\pi]$.  Then one has the equality:
\begin{equation}
\begin{aligned}
& \frac{1}{2 \pi} \int_0^{2\pi} d \theta |b(\psi)|^2 W^{k+1}_{\tau}(Y(S); E_*(z)) = \\
&\frac{1}{2 \pi} \int_0^{2\pi} d \theta \Sigma_{n,m} \exp(i
\theta(n-m))\langle\varphi_n, \varphi_m\rangle|
W^{k+1}_{\tau}(\widetilde{W}^0_0([-S,S]))
= \Sigma_n|\varphi_n|^2W^{k+1}_{\tau}.
\end{aligned}
\end{equation}
Combining with the inequality ($D^{-1}$ are bounded operators for all
$z \in \mathbb{C}$) $|b(\psi)|W^{k+1}_{\tau}(Y(S); E_*(z)) \leq C
|\widehat{\psi}|W^k_{\tau}$ for every $\theta$,
one gets the estimate:
$$\Sigma_n |\varphi_n|W^{k+1}_{\tau}(\widetilde{W}^0_0([-S,S])) 
   \leq C \Sigma_n |\psi_n|W^k_{\tau}(\widetilde{W}^0_0([-S,S])).$$
From the equality, $T^n(\psi) = \psi_n$, one has the desired estimate:
$$|R(\psi)|^2W^{k+1}_{\tau}(\widetilde{Y}(S); \widetilde{E}_*) = 
\Sigma_n |\varphi_n|^2 W^{k+1}_{\tau}(\widetilde{W}^0_0([-S,S])) \leq C\Sigma_n |D(\varphi_n)|^2 W^k_{\tau}$$
$$= C \Sigma_n |T^n(\psi)|^2 W^k_{\tau}(\widetilde{W}^0_0([-S,S])) = C |\psi|^2 W^k_{\tau}.$$
This completes the proof.
\end{proof}

\subsection{Computation of parametrized cohomology groups}
In 4.C, one will verify the assumption in proposition 4.1.
\begin{sublem}
$D_i\co W^{k+1}_{\tau}(Y(S); E_i(z)) \to W^k_{\tau}(Y(S); E_{i+1}(z))$ has 
closed range.
\end{sublem}
\begin{proof}
It is true for $z=1$. Recall $Y(S) = \widetilde{W}_0([-S,S])/ N_0(S)
\sim M_1(S)$.  Let $H_i \subset W^{k+1}_{\tau}(Y(S); E_i(z))$ be closed
subsets satisfying:

(1)\qua $H_1 \cup H_2 = W^{k+1}_{\tau}(Y(S); E_i(z))$ and

(2)\qua $\Supp H_1 \subset \widetilde{W}_0([-(S-1), S-1])$
and $\Supp H_2 \subset N_0([-S, -(S-1)]) \cup M_1([S-1,S])$.

By regarding $H_1 \subset W^{k+1}_{\tau}(Y(S); E_i(1))$, it follows that
$D_i(H_1)$ is a closed subspace of $W^k_{\tau}(Y(S); E_{i+1}(z))$. Let
$u \in H_2$.  Then one may associate $\widetilde{u} \in
W^{k+1}_{\tau}(Y(S); E_i(1))$ by $\widetilde{u}|N_0([-S, -(S-1)]) =u$
and $\widetilde{u}|M_1([S-1,S]) = z^{-1}u$.
Clearly $D_i(H_2) = D_i(\widetilde{H}_2)$ is a closed subspace of 
$W^k_{\tau}(Y(S); E_{i+1}(1))$.
Since this assignment is isometric, it follows that $H_2$ is also a
closed subspace of $W^k_{\tau}(Y(S); E_{i+1}(z))$.

This completes the proof.
\end{proof}
\begin{prop}
Let $\{ D_*\}$ be an AHS complex.  Suppose it gives a Fredholm complex and
all the cohomology groups vanish,
$H^i(W^k_{\tau}(Y(S); \Lambda^* ); D_*) =0$, $i=0,1,2$.
Then for $C(1) = \{ z: |z| =1 \}$, 
one has $H^*(W^k_{\tau}( Y(S) ; \Lambda^* ); (D_*)_z) =0$
for all $z \in C(1)$ and  $i=0,1,2$.
\end{prop}
\begin{proof}
The first part is essentially \cite[page 390]{taubes}.
Let $u \in W^k_{\tau}(Y(S); \Lambda^*_z )$ with $d_0(u) =0$.
Then passing through the identification 
$z^t\co \Lambda^0_z \cong  \Lambda^0$,
one has $d(z^{-t}u) =0$ where $z^{-t}u \in 
W^k_{\tau}(Y(S); \Lambda^0)$.
Clearly this shows $u=0$.

Next let us take $\alpha \in W^k_{\tau}(Y(S); \Lambda^1_z)$
with $d_+(\alpha) =0,  d^*_{\tau}(\alpha) =0$.
Then by integration by parts, one has $d(\alpha)=0$.
$\alpha$ is zero in $H^1(\widetilde{W}_0;\mathbb{R})$
(see the proof of  lemma 5.1).
Thus one may express $\alpha = d(f)$ on $\widetilde{W}_0$
for $f \in C^{\infty}(\widetilde{W}_0([-S,S]))$
and $d^*_{\tau} d(f) =0$.
Let $i\co N_0(S) \cong M_1(S)$ be the identification.
Then one has $i^*(f) = zf +\const$.
If one takes $f$ so that it vanishes at infinity, then $\const =0$.
Thus one has the equality:
$$0 = \int_{\widetilde{W}_0([-S,S])}\exp(\tau) 
  \langle f , d^*_{\tau} d(f)\rangle \vol = - \int \exp(\tau) |df|^2 .$$
This shows $df=0$.  Let $ O \subset C(1)$ be  the subset satisfying
a property that any $\alpha \in W^k_{\tau}(Y(S); \Lambda^1_z)$ with  
$d_+(\alpha) =0$ satisfies $d \alpha =0$.
The above implies $O$ is non-empty,   open and closed in $C(1)$.
Thus  $H^*(W^k_{\tau}( Y(S) ; \Lambda^* ); (D_*)_z) =0$, $*=0,1$
for all $z \in C(1)$. 

Now let us put $D= (D_0)^*_{\tau} \oplus D_1$, and consider $\Delta =
D^*_{\tau} D$ on $L^2_{\tau}(Y(S); \Lambda^1(z))$. Then there is a
constant $\lambda >0$ with $\lambda_z  > \lambda$, where $\lambda_z$
is the bottom of $\Delta$.  Suppose $H^2 \ne 0$ at $z = \exp(i \theta_0)$,
but $H^2 =0$ for all $z = \exp(i \theta)$, $0 \leq \theta < \theta_0$.
Then by choosing some $\theta$, it follows that there exists 
$u \in W^1_{\tau}(Y(S); \Lambda^2_z)$, $z= \exp(i \theta)$ with (1) $|u|W^1_{\tau} =1$ 
and (2) $|(D_1)^*_{\tau}(u)|L^2_{\tau} < \epsilon$.

Since $D_1$ is surjective, one finds $\alpha \in W^2_{\tau}(Y(S); \Lambda^1_z)$ with 
$D_1(\alpha) =u$ and $(D_0)^*_{\tau}(\alpha) =0$. Then there is a constant $C$
independent of $z$ such that:
$$|u|W^1_{\tau} \leq C |\alpha|W^2_{\tau} \leq C|(D_1)^*_{\tau}D_1(\alpha)|L^2_{\tau}.$$
This is a contradiction. This shows that the second cohomology also
vanishes for all $z \in C(1)$. This completes the proof. 
\end{proof}

\subsection{Computation of cohomology groups with $z=1$}
Now we compute the cohomology $H^i(W^k_{\tau}(Y(S); \Lambda^*); D_*) $ 
of the AHS complex for the simplest kinky handles $(W_0, N_0,M_1)$.
Notice that the end of $Y(S)$ is isometric to $ M \times [0, \infty)$
for a compact Riemannian three manifold $M$.
We will denote  $M_t = M \times \{t \} \subset Y(S)$.
\begin{prop}
Let $Y(S)$ be as   above.
Then for all $i=0,1,2$, the cohomologies
$H^i(W^k_{\tau}(Y(S); \Lambda^*); D_*)$ vanish. 
\end{prop}
\begin{proof}
Clearly $H^0(\AHS)=0$. Let us consider $H^1(\AHS)$. 
Let us take any representative 
$u $. By integration by parts, one finds $du=0$.
One may express $u= w + f\,dt$ around the end $Y(S)  \cong 
M \times[0, \infty)$, where $w$ does not contain $dt$ component.
 Then through this isometry, one considers $w_t \in \Lambda^1(M)$ for $t \in [0, \infty)$.
Let $d_3$ be the differential over $M$. Then for every $t$, one has $d_3w_t =0$.
Thus one may express $w_t=  \Sigma_i g_i(t) \alpha_i + d\mu_t$,
where $ \{ \alpha_i \}_i$ consists of the orthogonal  basis of 
$H^1(M;\mathbb{R})$. 
Here one may assume $d_3^*(\alpha_i) =0$.
By taking into account of $dt$ component of $du$,
one has the following  equation:
$$ \Sigma_i g_i(t)' \alpha_i + d_3\mu(t)' = d_3f_t.$$
This shows $g_i(t)' \equiv 0$.   Since $u \in L^2_{\tau}$, 
one concludes $g_i(t) \equiv 0$.
Then one has another equality, $f_t= \mu(t)' + c(t)$, 
where $c(t)$ are constants depending on $t$.
Thus one gets $u= d_3\mu_t + (\mu_t' +c(t))dt$.
\begin{sublem}
If $d_3\mu_t \in W^k_{\tau}$, then for a smooth family of constants $d(t)$,
one has $\mu_t -d(t) \in W^k_{\tau}$.
\end{sublem}
{\bf Proof of sublemma}\qua
Let us put $d(t) = (\vol_M)^{-1}\int_M \mu_t \vol$. Then 
we show that this family is the desired one.
Notice that one has the following bound: 
$|d_3\mu_t|L^2_{\tau} \geq C |\mu_t-d(t)|L^2_{\tau}$ for some positive constant $C$.
In particular one has  $\mu_t-d(t) \in L^2_{\tau}$.
Notice that one has the equality $\mu_t-d(t) =( \Delta_0)^{-1}d^* d(\mu_t-d(t))$,
and $( \Delta_0)^{-1}d^*$ is a translation invariant bounded operator.
Thus we get:
$$ \frac{\mu_t-d(t) }{dt} =
\frac{d}{dt} [ ( \Delta_0)^{-1}d_3^* d_3(\mu_t-d(t))]= ( \Delta_0)^{-1}d^*
\frac{d}{dt}[ d_3\mu_t].$$
 From the last term, one sees that $ \frac{\mu_t-d(t) }{dt}  \in L^2_{\tau}$.
By a similar consideration, one gets the result.
This completes the proof of the sublemma.
\endproof

{\bf Proof of proposition (continued)}\qua
Then by replacing $\mu_t$ by $\mu_t-d(t)$,  one may assume $\mu_t \in W^k_{\tau}$.
Thus $d\mu_t =d_3\mu_t+ \mu'_tdt \in W^{k-1}_{\tau}$.
On the other hand, let us consider the equality $f_t = \mu_t'+c(t)$.
Then it follows $c(t) \in W^{k-1}_{\tau}$. Let us put $C(t) = \int_0^t c(s)ds - \int_0^{\infty}c(t)$.
Then by \cite[lemma 5.2]{taubes} (see sublemma 3.1),
 one finds that $C(t) \in  W^k_{\tau}$, and $dC(t) = c(t)$.
In particular we have $u|\text{ end}Y(S) = df$, $f\in W^k_{\tau}$.
Then using a cut-off function, one may represent $u$ by a compactly supported 
smooth $1$ form, which is itself exact by a compactly supported smooth function,
since $H^1_{\text{cp}}(Y(S); \mathbb{R})=0$.
This shows $H^1(\AHS)=0$.

Next we consider $H^2(\AHS)$. Let us take a representative $u \in H^2(\AHS)$. One may choose $u$
so that it satisfies $d(e^{\tau}u)= e^{-\tau}d(ue^{\tau})=0$.
Since $H^2(Y(S); \mathbb{R})=0$, one may express $e^{\tau}u=d\mu$, $\mu \in C^{\infty}(\Lambda^1)$.
Let us denote $\mu|\operatorname{end}Y(S) =\beta + f\,dt$, where $\beta$
does not contain $dt$ component. Then we have the following relation:
$$d\mu|\operatorname{end}Y(S) = d_3\beta_t +(d_3f -\beta_t')dt =
d_3\beta_t + *_3d_3\beta_t \wedge dt$$
where $*_3$ is the star operator over $M$.
Let us decompose $\beta_t= \beta_t^1+ \beta_t^2$,
where $\beta_t^1 (\beta_t^2) $ does (not) consists of a closed form over $M$.
Then from the last two terms, one finds $d_3f_t =(\beta_t^1)'$.
In particular one may represent:
$$e^wu|\operatorname{end}Y(S) =d\beta_t^2 =d_3\beta_t^2 - (\beta_t^2)'
  \wedge dt = d_3\beta_t^2 + *_3d_3\beta_t^2 \wedge dt.$$
By the decomposition, one finds a positive constant $C$ such that:
$$|d_3\beta_t^2|W^{k-1}( M_t ) \geq C|\beta_t^2|W^k( M_t).$$
Now we have the next relations (put $\mu =\beta_t^2$ on the end):
$$e^{\tau}u =d\mu, \quad |\mu_t|W^k(M_t) \leq C|e^{\tau}u|W^{k-1}(M_t).$$
\begin{sublem} If $\delta >0$ is sufficiently small with respect to $S$,
then $e^{\tau}u \in L^2(Y(S); \Lambda_+^2)$.
\end{sublem}
{\bf Proof of sublemma}\qua
Let us consider $e^{\tau} u|\operatorname{end} =d \mu$.
One may assume that for every $t$,
 $\mu_t \in C^{\infty}(\Lambda^1(M_t))$
lies on the orthogonal complement of $\ker d$. Then $*_3 d$ is invertible
on $(\ker d)^{\perp}$ by corollary 3.1.
Moreover it is self-adjoint with respect to the $L^2$ inner product.
Since $\mu$ satisfies the equation, $(\frac{\partial}{\partial t} +
*_3d_3)\mu =0$, one has the exponential decay estimate for $\mu$.
More precisely there exist constants $C >0, \lambda_0 >0$ 
which are independent of $\delta$, such that:
$$|\mu|L^2(M_t) \leq \exp( -\lambda_0 t) \sup \{ |\mu|L^2(M_s);
  0 \leq s \leq 2t \}.$$
Notice that a priori, $\mu$ satisfies the following bound of
its growth
$|\mu|L^2(M_t) \leq C \exp(t \frac{\delta}{2})|u|L^2(M_t)$. 
Combining these estimates, one gets the exponential estimate for $\mu$.
Then one has also the exponential decay estimate for $d \mu$ on the end.
This completes the proof of the sublemma.
\endproof

{\bf Proof of proposition (continued)}\qua
Let $\varphi_t\co Y(S) \to [0, 1]$ be a cut-off function such that
$\varphi | M \times [t+1, \infty)  \equiv 0$, 
$\varphi|(M \times [t, \infty))^c \equiv 1$. Let us fix a large  $S$, and
 let $\delta >0$ satisfy the above condition.
By the above estimate on $W^k$ norm, 
one finds $\varphi_t \mu \in L^2 (Y(S); \Lambda_+^2) $ for all $t$.
Then for any small $\epsilon >0$, one finds a large $t$ so that:
$$\int |e^wu|^2 = \left|\int d(\varphi_t \mu ) \wedge d(\varphi_t \mu) +
  \int |e^wu|^2\right| < \epsilon.$$
This shows $u \equiv0$, and we have shown $H^2(\AHS)=0$ for this pair $(\delta,S)$.
This completes the proof of the proposition.
\end{proof}

\begin{rem}
{\rm
Suppose for all $\delta' \in [\delta, \delta_0]$, the differential of the
AHS complex has closed range, where $\delta$ is sufficiently small. Then
$H^*(\AHS)$ also vanishes for
all $\delta'$. This is seen as follows. Notice that for any choice
of the constants, one knows $H^0(\AHS) = H^1(\AHS) =0$. We want to calculate
$H^2(\AHS)$ when we vary $\delta' >0$. Let us take an isomorphism:
$$I\co W^*_{\tau(\delta)}(Y(S); \Lambda^*) \cong
  W^*_{\tau(\delta')}(Y(S); \Lambda^*)$$
where we denote $\tau(\delta)$ to express the weight constant. Passing
through $I$, one has $\{ W^*_{\tau(\delta)}(Y(S); \Lambda^*); I^{-1}
\circ d_* \circ I \}$, a continuous family  of Fredholm complexes. In
particular, the indices of these complexes are invariant.  Since for all cases,
one has $H^0(\AHS)= H^1(\AHS)=0$, one concludes $H^2(\AHS)=0$.}
\end{rem}

\begin{rem}
{\rm Suppose a cylindrical four--manifold $Y$ has nonzero
$H^2(Y;\mathbb{R})$.  Then in this case, one has a bound:
$$\dim H^2(\AHS) \leq 2 \dim H^2(Y;\mathbb{R}).$$
This directly follows from the proof of proposition 4.3.}
\end{rem}

\begin{cor}
Let us choose a sufficiently small $\delta >0$. Then 
one has an invertible operator:
$$P_w\co W^{k+1}_w(\widetilde{Y}(S); \Lambda^1) \cong
 W^k_w(\widetilde{Y}(S); \Lambda^0 \oplus \Lambda^2_+).$$
\end{cor}

\section{An asymptotic method to compute cohomology}
\subsection{$P_w$ over $Y(S,2)$}
Let $(W_0,N_0,M_1)$ be a simplest kinky handle, and $\widetilde{W}_0$ be
the Riemannian manifold constructed before.
In the previous section, one has the invertible operator:
$$P_{\tau}\co W^{k+1}_{\tau}( Y(S) ; \Lambda^1 ) \to 
W^k_{\tau} (Y(S); \Lambda^0 \oplus \Lambda_+^2 ).$$
 By the Fourier--Laplace transform, one gets an invertible
 Fredholm $P_{\tau}$ over $\widetilde{Y}(S)$ where:
$$\widetilde{Y}(S) =
\dots  \cup \widetilde{W}_0^{i-1}([-S,S]) \cup \widetilde{W}_0^i([-S,S]) \cup
\dots \quad i \in \mathbb{Z}.$$

Let $T_2$ be the periodic tree as $\mathbb{R} \cup_{n \in \mathbb{Z}} \mathbb{R}_0$,
where $\mathbb{R}_0$ is half the real line.
The  aim here is to show  that $P_{\tau}$ is invertible over
the periodic cover of $Y(S,2)$ defined below which corresponds to
$T_2$.
Let $(W_2, N_0,M_1,M_2)$ be a kinky handle with two kinks, and
$\widetilde{W}_2$ be the corresponding complete Riemannian manifold.
Thus one has three ends in the horizontal direction,
$N_0([0, \infty)), M_j([0, \infty))$, $j=1,2$.
Let us choose a large $S$, and put:
\begin{equation}
\begin{aligned}
& \widetilde{Y}(S)_0 = \widetilde{W}^1_0((- \infty, S]) \cup \widetilde{W}^2_0([-S,S]) \cup \dots, \\
& Y(S,2)' =   \widetilde{W}_2 \backslash N_0([S, \infty)) \cup M_1([S, \infty))
/ N_0(S) \sim M_1(S), \\
& Y(S,2) =   Y(S,2)' \backslash M_2([S, \infty)) \cup \widetilde{Y}(S)_0 \backslash N_0([S, \infty)).
\end{aligned}
\end{equation}
One may equip  a weight function $\tau$ on $Y(S,2)$ as before.
\begin{lem}
$d_*\co W^{k+1}_{\tau} ( Y(S,2) ; \Lambda^* ) \to 
W^k_{\tau}( Y(S,2) ; \Lambda^{*+1})$ 
has closed range with $H^0(\AHS)=H^1(\AHS)=0$.
\end{lem}
\begin{proof}
Closedness follows from the excision method used before.
Let us take $(W_s= S^1 \times D^3, N_0 =S^1 \times D^2)$, where $N_0$ represents a generator
of $\pi_1(W_s)$. By the construction in 2.A, one can get a complete
Riemannian manifold $\widetilde{W}_s$
with one end $N_0([0, \infty))$ along the horizontal direction. Let us put:
$$Y(S,2)(0) = Y(S,2)' \backslash M_2([S, \infty)) \cup \widetilde{W}_s \backslash N_0([S, \infty)).$$
Notice that $Y(S,2)(0)$ is diffeomorphic to $Y(S, 2)'$,
but Riemannian metrics are different. ($Y(S,2)(0)$ has cylindrical
end.) By corollary 3.2, the differential of the AHS complex has closed range.

Let $H_i \subset W^{k+1}_{\tau}(Y(S,2); \Lambda^*)$ be closed subsets with:

(1)\qua $H_1 \cup H_2 = W^{k+1}_{\tau}(Y(S,2); \Lambda^*)$,

(2)\qua $\Supp H_1 \subset Y(S,2)' \backslash M_2([S, \infty))$ and
$\Supp H_2 \subset \widetilde{Y}(S)_0 \backslash N_0([S-1, \infty))$.

One regards $H_1 \subset W^{k+1}_{\tau}(Y(S,2)(0); \Lambda^*)$
and $H_2 \subset W^{k+1}_{\tau}(\widetilde{Y}(S); \Lambda^*)$.
From this, it follows that $d_*(H_i)$ are both closed subspace of
$W^k_{\tau}(Y(S,2); \Lambda^{*+1})$.

Let us see $\ker P_w=0$. 
$\ker P_w$ is isomorphic to the first cohomology $H^1(\AHS)$.
One can easily check  $H^1(Y(S,2); \mathbb{R}) = \mathbb{R}$.
Let  us take a nonzero  element $w $. Then this has the property 
$\langle w, C\rangle \ne 0$ for any loop $C $ representing a generator
of $\pi_1(Y(S,2))$.  One may choose $C$ sufficiently near to infinity
while having a bounded length.

Let us take $u \in H^1(\AHS)$, and consider its class $[u] \in 
H^1(Y(S,2); \mathbb{R})$. It follows $[u] =0 $ , since one can make
$\langle u, C\rangle$ arbitrarily small by choosing $C$ as above.
Thus one may express $u = df$ for $f \in C^{\infty}(Y(S,2))$.
One may assume $f \in L^2_w(Y(S,2))$, by subtracting some constant.
This shows $u =0$.  Thus we have shown the result.
\end{proof}

\subsection{Some continuity of  $H^1({\rm AHS})$}
Let us introduce another weight function:
$$\mu\co \widetilde{Y}(S) \to [0, \infty)$$
by $\mu|\widetilde{W}_n(s \leq S)(x) = |n| + t(x)$.
We choose another small constant $\delta'$ (we take it so that $\delta'
\ll \delta$,
where $\delta $ is the weight constant for $\tau$).
Then we introduce another weight $w = \delta' \mu + \tau$ and weighted 
Sobolev $k$ spaces $W^k_w$, where:
$$(|u|W^k_w)^2 = \int_{\widetilde{Y}} \exp(\delta' \mu +
  \tau)(\Sigma_{l=0}^k |\nabla^l u|^2).$$
Since one knows $H^*(\AHS) =0$ over $\widetilde{Y}(S)$, it follows that
there exists a positive constant $C>0$ such that $|P_w|, |(P_w)^*_w| \geq C$, for any 
small $\delta' \in [0, \delta_0']$.

Let $Y$ be a complete Riemannian manifold of bounded geometry.
Let us take a family of smooth maps:
$$w(\delta')\co Y \to [0, \infty), \quad w(\delta'') (x) \leq w(\delta')(x),
\quad \delta'' \leq \delta' \in [0, \delta], \quad x \in Y.$$
Suppose that for all $w(\delta')$, the corresponding AHS complex between
weighted Sobolev spaces are Fredholm. In any case it is immediate to see	
$H^0(\AHS)=0$. 
Here one has some continuity property as:
\begin{lem}
$H^1(\AHS) =0$ for $\delta' =0$ when $H^1(\AHS) =0$ for $\delta' >0$.
\end{lem}
\begin{proof}
Let us denote $\ker (d_+)^*_{w( \delta')} = H^2(\delta')$ 
 (the second cohomology $H^2(\AHS)$ when it is defined).
Let us choose sufficiently small constants:
$$\delta \gg \delta' > \delta'' \geq 0.$$
Then one has $H^2( \delta'') \cong H^2( \delta') $.
This follows from that by varying these weights from $w( \delta')$
to $w( \delta'')$, one gets a family of Fredholm complexes between
the weighted Sobolev spaces 
(these spaces also vary with respect to the deformation of the weights).
For every value of weights, one has $H^0 (\AHS) =H^1(\AHS)=0$.
Thus one gets the above statement.
 
\begin{sublem}
If $H^2(\delta'') \ne 0$, then one has also 
$H^2( \delta') \ne 0$.
\end{sublem}
{\bf Proof of sublemma}\qua
Let us denote $w= w(\delta') $ and 
$w' = w( \delta'') $.
Let us take $v \in L^2_{w'} $ with $d^*(\exp(w')v) =0$.
Then one puts $u = \exp(-w+w')v$. One may assume $u \in L^2_w$.
Since $u$ satisfies $d^*(\exp(w)u)=0$, one gets the result.
\endproof
{\bf Proof of lemma (continued)}\qua
Recall that by putting $w = w(\delta')$,
one also has a Fredholm complex
$\{ W^*_w(Y(S,2); \Lambda^*), d^* \}$.
Then by varying a parameter $\delta' \in [0,\delta]$,
one has a family of Fredholm complexes:
 $$\{ W^*_{w(\delta')}(Y(S,2); \Lambda^*), d^* \}.$$
 From the above proof, one has the inclusion $H^2(\delta'')  \subset H^2(\delta')$
for all $0 \leq  \delta'' \leq \delta' $.
Now let us see $H^2(\delta') = H^2(0)$ and $H^1(0) =0$.
Suppose $H^1(0) \ne 0$. Then by invariance of Fredholm indices, 
one must have $\dim H^2( 0) > \dim H^2(\delta')$. This contradicts
the above.  This completes the proof.
\end{proof}

\subsection{Computation of $H^2({\rm AHS})$}
To show that  $P_{\tau}\co W^{k+1}_{\tau}(Y(S,2); \Lambda^1) \to
W^k_{\tau}(Y(S,2); \Lambda^0 \oplus \Lambda^2_+)$ is invertible,
one uses an asymptotic method.
Roughly speaking one approximates $Y(S,2)$ by a family of Riemannian manifolds with cylindrical ends.
Then the spectrum of $P_{\tau}$ over each cylindrical manifold
has a uniform lower bound on $1$ forms. From this one gets a lower bound
of $(d_+)^*_{\tau}$ over $Y(S,2)$. In the presence of $H^2(Y;\mathbb{R})$, 
one can get a uniform bound of $\dim H^2(\AHS)$. The approximation of
spaces corresponds to the one of infinite tree by its finite
subtrees. This method completely works for higher stages in 5.D.

Computation of $\dim H^2(\AHS)$ uses information of $H^1(\AHS)$ on both
approximation spaces and the limit space. The former is obtained by
proposition 4.3 and the latter by sublemma 3.1. On the other hand for
all cases $H^0(\AHS)=0$.  Then using this and proposition 4.3, one can
apply the asymptotic method in this section to verify $H^1(\AHS)=0$
for the limit space. This will be another method to compute $H^1(\AHS)$
without using sublemma 3.1.

Let us take $(W_s = S^1 \times D^3, N_0 = S^1 \times D^2)$
where $N_0$ represents a generator of $\pi_1(W_s)$.
Then  by the construction in 2.A, one gets a complete Riemannian
manifold $\widetilde{W}_s$ with one end $N_0([0, \infty))$ along
horizontal direction.
Let $(W_0,N_0,M_1)$ be the simplest kinky handle as before.
Let us put:  
\begin{equation}
\begin{aligned}
& Y(S)_0(n)' = \widetilde{W}_0^1([-S,S]) \cup_{M_1^1(S) \sim N_0^2(S)} \widetilde{W}_0^2([-S,S]) \\
& \qquad         \dots \cup_{M_1^{n-1}(S) \sim N_0^n(S)} \widetilde{W}_0^n([-S,S]), \\
& Y(S)_0(n) = Y(S)_0(n)' \cup_{M_1^n(S) \sim N_0(S)} \widetilde{W}_s \backslash  N_0([S, \infty)), \\
& Y(S)_0'' = \dots  \widetilde{W}_0^{-n}([-S,S]) \cup_{M_1^{-n} (S) \sim N_0^{-n+1}(S)} \widetilde{W}_0^{-n+1}([-S,S])  \dots \\
& \qquad  \cup_{M_1^{n-1}(S) \sim N_0^n (S)} 
   \widetilde{W}_0^n([-S,S]) 
           \cup_{M_1^n(S) \sim N_0(S)} \widetilde{W}_s \backslash  N_0([S, \infty)).
\end{aligned}
\end{equation}
Notice $Y(S)_0(\infty)' = Y(S)_0$.
Let us put $Y(S,2)(n)' = Y(S,2)' \cup_{M_2(S) \sim N_0^0(S)} Y(S)_0(n)'$ and:
$$Y(S,2)(n) = Y(S,2)' \cup_{M_2(S) \sim N_0^0(S)} Y(S)_0(n).$$
$Y(S,2)(n)$ is a complete Riemannian manifold with cylindrical end.
Let the weight function be $\tau\co Y(S,2)(n) \to [0, \infty)$, with weight 
$\delta$
as before. 
\begin{prop}
One gets isomorphisms:
$$P_{\tau}\co W^{k+1}_{\tau} (Y(S,2)(n); \Lambda^1) \cong
  W^k_{\tau}(Y(S,2)(n); \Lambda^0 \oplus \Lambda^2_+).$$
\end{prop}
\begin{proof}
By proposition 4.3, the conclusion is true when one uses a sufficiently
small weight constant $\delta'$. On the other hand, by the first part of the proof
of lemma 5.1, the differential of the AHS complex has closed range for all
$\delta' \in (0, \delta]$. From this, one gets $H^0(\AHS) = H^1(\AHS) =0$
for all $\delta'$.  To see $H^2(\AHS)=0$, one can use the same method as
proposition 4.2.  This completes the proof.
\end{proof}

In particular there are constants $C_n$ such that 
for $\Delta = (P_{\tau})^*_{\tau} \circ P_{\tau}$, 
one has the bounds
 $|u|W^{k+2}_{\tau}(Y(S,2)(n)) \leq C_n|\Delta(u)|W^k_{\tau}(Y(S,2)(n)) $.
Also one has  bounds: 
$$ |u|W^{k+2}_{\tau}(Y(S,2)) \leq C|\Delta(u)|W^k_{\tau}(Y(S,2)),$$
$$ |u|W^{k+2}_{\tau}(Y(S)_0'') \leq C |\Delta(u)|W^k_{\tau}(Y(S)_0'').$$
\begin{sublem}
There is a lower bound $C_n \leq C$  for all $n$.
\end{sublem}
\begin{proof}
For any small $\epsilon >0$, there exists a large $n$
with the following property:
for any $u \in W^{k+2}_{\tau}(Y(S,2)(n))$ with $|u|W^{k+2}_{\tau}  =1$,
there is  at least one $n'$, $0 \leq n' \leq n$ with
$|u|W^{k+2}_{\tau}(\widetilde{W}_0^{n'}([-S,S]) < \epsilon$.
Let $\varphi(n')$ be a cut-off function on $Y(S,2)(n)$
with $\varphi(n')| ( \widetilde{W}_0^{n'}([-S,S]))^c \equiv 1$
and $Y(S,2)(n) \backslash \Supp d\varphi(n')$ has two components.
 Then one may express $\varphi(n')u= u_1+ u_2$.
One may regard $u_1 \in W^{k+2}_{\tau}(Y(S,2))$, 
 $u_2 \in W^{k+2}_{\tau}(Y(S)_0'')$.
In particular one has the estimates: 
\begin{equation}
\begin{aligned}
 |u|W^{k+2}_{\tau} & \leq |u_1|W^{k+2}_{\tau} + |u_2|W^{k+2}_{\tau}
+ |(1- \varphi(n'))u|W^{k+2}_{\tau} \\
& \leq C\{ |\Delta(u_1)|W^k_{\tau} + |\Delta(u_2)|W^k_{\tau} \} + \epsilon 
 \leq C|\Delta(u)|W^k_{\tau} + C \epsilon
\end{aligned}
\end{equation}
where $C$ is independent of $\epsilon$ and $u$.
This gives the result.
\end{proof}

Now suppose $H^2(\AHS)$ is nonzero over $Y(S,2)$, 
and take $u \in \ker(d_+)^*_{\tau} \cap W^{k+1}_{\tau}(Y(S,2);
\Lambda^2_+)$ with $|u|W^{k+1}_{\tau} =1$.
For each $n$, let us take a cut-off function $\varphi(n)$ on $Y(S,2)$ 
with $\varphi(n)| Y(S,2)(n-1)' \equiv 1$, $\varphi(n)| (Y(S,2)(n)')^c \equiv 0$.
Then one may regard $\varphi(n)u \in W^{k+1}_{\tau}(Y(S,2)(n); \Lambda^2_+)$.
Thus there exists $v_n \in W^{k+2}_{\tau}(Y(S,2)(n); \Lambda^1)$
with $d^*_{\tau}(v_n) =0$, $d_+(v_n)= \varphi(n)u$.
By the above sublemma, one has the uniform  estimates: 
\begin{equation}
\begin{aligned}
&    1- \epsilon \leq |\varphi(n)u|W^{k+1}_{\tau}   = |d_+(v_n)|W^{k+1}_{\tau} \\
& \leq C  |v_n|W^{k+2}_{\tau} \leq C| d_+(v_n)|W^{k+1}_{\tau}
\leq  C|\varphi(n)u|W^{k+1}_{\tau}. \quad (*)
\end{aligned}
\end{equation}

By the above $(*)$, one has the uniform estimates 
$ C \leq |v_n|W^{k+2}_{\tau} \leq C'$. 
Let us take a sequence $1 \gg \epsilon_1 \gg \epsilon_2 \gg \dots \to 0$.
Then for each $n$, there exists a large $N = N(n)$ and  $L(n) \gg 0$
such that at least one $|v_N|W^{k+2}_{\tau}(\widetilde{W}_0^{n'}([-S,S])) $
is less than $\epsilon_n$ for $N- L(n) \leq n' \leq N$.
Let us take a subsequence $\{ v_{N(n)}\}_n$
and denote it by $\{ v_n \}_n$. 
For simplicity of the notation, one may assume $n' = n-1$.

Using this, one can verify $\ker(d_+)^*_{\tau} =0$ over $Y(S,2)$.
Let us put $\Delta = d d^*_{\tau} \oplus (d_+)^*_{\tau} d_+$
on $W^{k+2}_{\tau}(Y(S,2); \Lambda^1)$.
 From above, one has the estimate:
\begin{equation}
\begin{aligned}
  |dd^*_{\tau} &  (\varphi  (n-1)v_n)|W^k_{\tau} \\
& \leq  C \{   |dv_n|W^k_{\tau}(\widetilde{W}_0^{n-1}([-S,S]) ) 
 +|v_n|W^k_{\tau}(\widetilde{W}_0^{n-1}([-S,S])) \} 
\leq C \epsilon_n
\end{aligned}
\end{equation}
where $C$ is independent of $n$.
Then  one has the estimates:
\begin{equation}
\begin{aligned}
| \Delta & (\varphi(n-1)v_n)|W^k_{\tau}  \leq 
C[ |  (d_+)^*_{\tau}    ( \varphi(n-1)  \varphi(n) u)|W^k_{\tau} \\
& +    |    (d_+)^*_{\tau}  [d_+, \varphi(n-1) ]v_n|W^k_{\tau}(\Supp d \varphi(n-1))]   + C \epsilon_n\\
&  \leq C[ |\varphi(n-1)\varphi(n)  (d_+)^*_{\tau}(u)|W^k_{\tau} \\
      & \qquad      +  |[ (d_+)^*_{\tau} , \varphi(n-1)\varphi(n) ] u|W^k_{\tau} (\cup_{j=n-1,n} \widetilde{W}_0^j([-S,S]))  \\
&  \qquad +  |    (d_+)^*_{\tau}  [d_+, \varphi(n-1)]v_n|W^k_{\tau}(\widetilde{W}_0^{n-1}([-S,S]))
+ C \epsilon_n \\
& \leq  |[ (d_+)^*_{\tau} , \varphi(n-1)\varphi(n) ] u|W^k_{\tau} +
       |    (d_+)^*_{\tau}  [d_+, \varphi(n-1)]v_n|W^k_{\tau}
+ C \epsilon_n.
\end{aligned}
\end{equation}
One may assume that the last term is arbitrarily small.
This shows that there is a sequence $\{ w_n \in W^{k+2}_{\tau}(Y(S,2);
\Lambda^1) \}$ with $C \leq |w_n|W^{k+2}_{\tau} \leq C'$ and
$|P(w_n)|W^{k+1}_{\tau} $ converges to zero. This is a contradiction. This
shows $\ker(d_+)^*_{\tau} =0$ over $Y(S,2)$.

Let $(W_t, N_0) = (D^4, S^1 \times D^2)$ be the standard disk.
Let $Y = \widetilde{W}_t \backslash N_0([S, \infty)) \cup \widetilde{Y}(S)_0 \backslash N_0([S, \infty))$. In this case,
one has $H^0(\AHS) = H^1(\AHS)=0$ (notice $Y$ is simply connected). But $H^2(Y;\mathbb{R})=\mathbb{R}$.
In this case, one has the following:
\begin{cor} Suppose $Y$ has nonzero $H^2(Y;\mathbb{R})$. Then 
one has a bound $\dim H^2(\AHS) \leq 2 \dim H^2(Y;\mathbb{R})$.
\end{cor}
\begin{proof}
Suppose $\dim H^2(\AHS) \geq 2 \dim H^2(Y;\mathbb{R}) +1= b_2+1$.
Let us take $L^2_{\tau}$ orthogonal vectors 
$u_1, \dots, u_{b_2+1}, \dots  \in L^2_{\tau}(Y: \Lambda^2_+) \cap \ker
(d_+)^*_{\tau}$ with $|u_i|L^2_{\tau} =1$.

One can make a family of cylindrical manifolds $Y(n)$ by the above method.
By  a straightforward calculation, one has an upper bound of
dim$H^2(Y(n);\mathbb{R})$ by $b_2$.  By remark 4.2,
it follows that $\dim H^2(\AHS)$ has also an upper bound
over $Y(n)$. Now let $\varphi(n)$ be as above. Then for any small $\epsilon >0$,
there exists a large $n_0$ such that for all $n \geq n_0$, one has:
\begin{equation}
\begin{aligned}
& |(d_+)^*_{\tau}(\varphi(n)u_i)|L^2_{\tau} \leq \epsilon, \quad
|\varphi(n)u_i|L^2_{\tau} \geq 1- \epsilon, \\
& \langle \varphi(n)u_i, \varphi(n)u_j\rangle|L^2_{\tau} \leq \epsilon, \quad
i= 1, \dots, b_2+1.
\end{aligned}
\end{equation}
Let $v_1, \dots, v_l$ be an orthonormal basis of $H^2(\AHS)$ over $Y(n)$,
$l \leq b_2$. Then one may express $\varphi(n)u_i = \Sigma_j a^i_j v_j +d_+(\alpha_i)$,
$a^i_j \in \mathbb{R}$. Let $\overline{a}_i = (a^i_1, \dots, a^i_l) \in \mathbb{R}^l$.
Then the set of vectors $\overline{a}_1, \dots, \overline{a}_{b_2+1}$
would satisfy $|\overline{a}_i| \geq 1- \epsilon$ and
$\langle\overline{a}_i, \overline{a}_{i'}\rangle \leq \epsilon$
with respect to the standard norm in $\mathbb{R}^l$. Since $l < b_2+1$,
it is impossible to find such set. This completes the proof.
\end{proof}

\subsection{Fourier--Laplace transforms on higher stages}
Let us put:
\begin{equation}
\begin{aligned}
& \widetilde{W}_2([-S,S]) =   \widetilde{W}_2 \backslash N_0([S, \infty)) \cup_{j=1,2} M_j([S, \infty)), \\
& Y(S,2)_0 =   \widetilde{W}_2 \cup_{M_2(S) \sim N_0(S)}  \widetilde{Y}(S)_0 \backslash N_0([S, \infty)), \\
& \widetilde{Y}(S,2) =   \dots Y(S,2)_0^{-n} \cup_{M_1^{-n}(S) \sim N_0^{-n+1}(S)}
Y(S,2)_0^{-n+1} \dots \\
& \widetilde{Y}(S,2)_0 =  Y(S,2)_0^1  \cup_{M_1^1(S) \sim N_0^2(S)}
Y(S,2)_0^2 \cup_{M_1^2(S) \sim N_0^3(S)} Y(S,2)_0^3  \dots 
\end{aligned}
\end{equation}
where $(Y(S,2)_0^n , M_j^n(S), N^n_0(S))$ is the same triple $(Y(S,2)_0, M_j(S), N_0(S))$
as before.
Notice that $\widetilde{Y}(S,2)_0$ is diffeomorphic to $CH(T^0_2)$, where  
$T_2^0 = \mathbb{R}_+ \cup_{n  \in \mathbb{N}} \mathbb{R}_+$.
Now by the Fourier--Laplace transform, one gets the following:
\begin{cor}
$P_{\tau}\co W^{k+1}_{\tau}(\widetilde{Y}(S,2) ; \Lambda^1) 
\to W^k_{\tau}(\widetilde{Y}(S,2); \Lambda^0 \oplus \Lambda^2_+)$
is invertible.
\end{cor}

Now one defines:
$T_3 = \mathbb{R} \cup_{n \in \mathbb{Z}} T_2^0$ and
 $T_3^0 = \mathbb{R}_+ \cup_{n \in \mathbb{N}} T_2^0$. Similarly
 $T_4= \mathbb{R} \cup_{n \in \mathbb{Z}} T_3^0$. One inductively defines
$T_j$, $j =1,2, \dots$  as:
$$T_{j+1} =  \mathbb{R} \cup_{n \in \mathbb{Z}} T_j^0, \quad    
T_{j+1}^0 = \mathbb{R}_+ \cup_{n \in \mathbb{N}} T_j^0    $$
where one puts $T_1 = \mathbb{R}$, $T_1^0 = \mathbb{R}_+$.

One can construct the corresponding spaces.
Let $(W_2, N_0,M_1,M_2)$ be a kinky handle with two kinks.
One has already defined $Y(S)$, $Y(S,2)$. 
Let us define inductively  $Y(S,j)$ as follows:
\begin{equation}
\begin{aligned}
& Y(S,j)_0 =   \widetilde{W}_2([-S,S]) \cup_{M_2(S) \sim N_0(S)} \widetilde{Y}(S,j-1)_0, \\
& Y(S,j) =   Y(S,j)_0 / N_0(S) \sim M_1(S), \\
& \widetilde{Y}(S,j) =   \dots Y(S,j)_0^{-n} \cup_{M_1^{-n}(S) \sim N_0^{-n+1}(S)}
Y(S,j)_0^{-n+1} \dots \\
& \widetilde{Y}(S,j)_0 =   Y(S,j)^1  \cup_{M_1^1(S) \sim N_0^2(S)}
Y(S,j)^2 \cup_{M_1^2(S) \sim N_0^3(S)} Y(S,j)^3  \dots 
\end{aligned}
\end{equation}
One may express $CH(T^0_j)= \widetilde{Y}(S,j)_0$.
The previous method works for all $Y(S,j)$ iteratively.
\begin{prop} $P_{\tau}$ over $Y(S,j)$ are all Fredholm with
$H^*(\AHS) =0$ for $* =0,1,2$ and $j = 1,2, \dots$.
\end{prop}

In the notation in 2.B, one expresses $T^0_j = (T_{2,2, \dots , 2,1})_0$ 
($ (2, \dots,2)$ $j-1$ times).
Let $n_1, \dots, n_k \in \{ 1,2, \dots \}$ be
a set of positive integers. 
Then  using kinky handles with $n_j$ kinks, one has  a natural extension, and gets 
 $(T_{n_1, \dots, n_k,1})_0 $ which is a signed infinite tree. 
For Riemannian metrics on  $CH((T_{n_1, \dots,n_k, 1})_0)$,
see 2.B. 

Let $\overline{n} = \{n_1, \dots, n_{l-1},1\}$ be a set of positive
integers, and denote the corresponding
homogeneous tree of bounded type by $(T_{\overline{n}})_0$.
By the previous method, one gets a complete Riemannian metric and a
weight function on every $CH((T_{\overline{n}})_0)$.
Recall that  one has constructed
complete Riemannian metrics and weight functions on $Y(S, \overline{n})$.
\begin{prop}
$P_{\tau}\co W^{k+1}_{\tau}(Y(S, \overline{n}); \Lambda^1) \to 
W^k_{\tau}(Y(S, \overline{n}); \Lambda^0 \oplus \Lambda^2_+)$
gives an isomorphism for any  $\overline{n} = \{n_1, \dots, n_{l-1},1 \}$.
\end{prop}
\begin{proof}
One has shown the result for $l=2$. Suppose the result is true for $l \leq l_0$.
 Let us put $\overline{n} = \{ n_1, \dots, n_{l_0 +1},1\}$.  Then 
by Fourier--Laplace transform and excision method used  before,
 one knows  $P_{\tau}\co W^{k+1}_{\tau}(Y(S, \overline{n}); \Lambda^1) \to
  W^k_{\tau}(Y(S, \overline{n}); \Lambda^0 \oplus \Lambda^2_+)$
gives a closed operator with $H^*(\AHS) =0$ for $*=0,1$.
One may follow the same process to see $H^2(\AHS)=0$ as
5.C.
Thus one has shown the result for $l_0+1$. This completes the induction step.
\end{proof}

In practical applications,  one considers open four--manifolds composed of
one $0$--handle attached with Casson handles.
Recall that $k(S^2 \times S^2) \backslash pt$ is homotopy-equivalent to
some wedges of $S^2$. In particular  $H^1_{\text{cp}}(M;\mathbb{R}) =
H^3(M;\mathbb{R}) = 0$.  Recall also that 
it  has a link picture by $k$ disjoint union of Hopf links with $0$--framings.  Let    $(W_t, M_1, \dots, M_{2k}) = (D^4, S^1 \times D^2,
\dots, S^1 \times D^2)$ express the link diagram of  $k(S^2 \times S^2)
\backslash pt$.

Let $(T_1)_0, \dots, (T_{2k})_0$ be signed homogeneous trees of bounded type.
Let us consider an open four--manifold $S$ obtained by attaching $CH((T_i)_0)$ along
$(D^4, S^1 \times D^3, \dots, S^1 \times D^3)$.
Then by the previous procedure, one can equip a complete Riemannian metric on $S$:
 $$S \equiv   \widetilde{W}_t \natural \cup_l CH((T_l)_0).$$
As before one can also equip a weight function $\tau$, and  the AHS complex over $S$.
\begin{cor} The differential of AHS complex has closed range over $S$
with $H^0(\AHS)=H^1(\AHS) =0$ and $\dim H^2(\AHS) \leq 2 \dim
H^2(S;\mathbb{R})$.
\end{cor}

By the work of Freedman \cite{freedman},
the end of $S$ admits a topological color, $\cong S^3 \times [0, \infty)$.
In fact $S$ is homeomorphic to $k(S^2 \times S^2) \backslash pt$. 
Now we have completed  the verification that any open four--manifold with
a tree-like end of bounded type can admit an admissible pair $(g, \tau)$ on it.

\end{document}

%% file: gtoutput.tex
\def\ifplaintex{\expandafter\ifx\csname documentclass\endcsname\relax}

\ifplaintex 
\else
\fi

\expandafter\ifx\csname beginpicture\endcsname\relax
\expandafter\ifx\csname documentclass\endcsname\relax
\input pictex \else
\input prepictex \input pictex \input postpictex \fi\fi

\def\gt{{\mathsurround=0pt\it $\cal G\mskip-2mu$eometry \&\ 
$\cal T\!\!$opology}}        

\def\gtp{{\mathsurround=0pt\it $\cal G\mskip-2mu$eometry \&\ 
$\cal T\!\!$opology $\cal P\!$ublications}}  

\def\lognumber#1{\def\thelognumber{#1}}
\def\volumenumber#1{\def\thevolumenumber{#1}}
\def\papernumber#1{\def\thepapernumber{#1}}
\def\volumeyear#1{\def\thevolumeyear{#1}}

\def\pagenumbers#1#2{\def\startpage{#1}\def\finishpage{#2}}
\def\published#1{\def\publishdate{#1}}
\def\proposed#1{\def\theproposer{#1}}
\def\seconded#1{\def\theseconders{#1}}
\def\received#1{\def\receiveddate{#1}}
\def\revised#1{\def\reviseddate{#1}}
\def\accepted#1{\def\accepteddate{#1}}
\def\asciititle#1{\def\theasciititle{#1}}

\long\def\asciiabstract#1{\long\def\theasciiabstract{#1}}

\let\\\par\let\thelognumber\relax
\let\thevolumenumber\relax\let\thepapernumber\relax
\let\thevolumeyear\relax\let\thesamplenumber\relax\let\startpage\relax
\let\finishpage\relax\let\publishdate\relax\let\receiveddate\relax
\let\reviseddate\relax\let\accepteddate\relax\let\theasciititle\relax
\let\theasciiauthors\relax
\let\theasciiabstract\relax
\let\theasciiemail\relax\let\theshortauthors\relax\let\theshorttitle\relax

\long\def\maketitlep{   

\count0=\startpage

\gt\hfill      
\beginpicture
\setcoordinatesystem units <0.33truein, 0.33truein> point at 2.2 0.9
\setplotsymbol ({$\cal G$})
\plotsymbolspacing=9truept
\circulararc 315 degrees from 0 1 center at 0 0
\setplotsymbol ({$\cal T$})
\circulararc 315 degrees from 1 -1 center at 1 0
\endpicture
\break
{\small\ifx\thesamplenumber\relax 
Volume \else Sample
\fi\thevolumenumber\ (\thevolumeyear)
\startpage--\finishpage\nl
Published: \publishdate}
\vglue 0.5truein plus 0.4fil minus 0.1truein

{\parskip=0pt\leftskip 0pt plus 1fil\def\\{\par\smallskip}{\ifplaintex\large
\else\Large\fi\bf\thetitle}\par\medskip}   

\vglue 0pt plus 0.1fil 

{\parskip=0pt\leftskip 0pt plus 1fil\def\\{\par}{\sc\theauthors}
\par\medskip}

\vglue 0pt plus 0.1fil 

{\small\parskip=0pt\let\newline\\
{\leftskip 0pt plus 1fil\def\\{\par}{\sl\theaddress}\par}
\expandafter\ifx\theemail\relax    
\relax\else\vglue 5pt plus 0.02fil minus 2pt\def\\{\stdspace{\rm 
and}\stdspace} 
\cl{Email:\stdspace\tt\theemail}\fi
\ifx\theurl\relax                  
\relax\else\vglue 5pt plus 0.02fil minus 2pt\def\\{\stdspace{\rm 
and}\stdspace}
\cl{URL:\stdspace\tt\theurl}\fi\par}

\vglue 7pt plus 0.3fil minus 3pt

{\bf Abstract}
\vglue 5pt plus 0.1fil minus 2pt

\theabstract

\vglue 7pt plus 0.3fil minus 3pt

{\bf AMS Classification numbers}\quad Primary:\quad \theprimaryclass

Secondary:\quad \thesecondaryclass

\vglue 5pt plus 0.3fil minus 2pt

{\bf Keywords:}\quad \thekeywords

\vglue 10pt plus 0.5fil minus 5pt

{\small  Proposed: \theproposer\hfill Received: \receiveddate\nl
Seconded: \theseconders\hfill 
\ifx\reviseddate\relax                         
Accepted: \accepteddate                        
\else
Revised: \reviseddate                          
\fi}
\eject
}       

\let\maketitlepage\maketitlep
\let\maketitle\maketitlepage

\font\phead=cmsl9 scaled 950
\font=cmsl9 scaled 1050
\font\pnum=cmbx10 scaled 913
\font\lnum=cmbx10 
\font\pfoot=cmsl9 scaled 950
\font=cmsl9 scaled 1050
\ifplaintex
\headline{\vbox to 0pt{\vskip -4.5mm\line{\small\phead\ifnum
\count0=\startpage ISSN 1364-0380 (on line)
1465-3060 (printed) \hfill {\pnum\folio}\else\ifodd\count0\def\\{ }%
\ifx\theshorttitle\relax\thetitle\else\theshorttitle\fi\hfill{\pnum\folio}
\else\def\\{ and }{\pnum\folio}\hfill\ifx\theshortauthors\relax\theauthors
\else\theshortauthors\fi\fi\fi}\vss}}
\footline{\vbox to 0pt{\vglue 0mm\line{\small\pfoot\ifnum\count0=\startpage
\copyright\ \gtp\hfill\else
\gt, Volume \thevolumenumber\ (\thevolumeyear)\hfill\fi}\vss
}}
\else
\makeatletter
\def\@oddhead{{\small\lhead\ifnum\count0=\startpage ISSN 1364-0380 (on line)
1465-3060 (printed) \hfill {\lnum\number\count0}\else\ifodd\count0
\def\\{ }\ifx\theshorttitle\relax \thetitle \else\theshorttitle\fi\hfill
{\lnum\number\count0}\else\def\\{ and }{\lnum\number\count0}
\hfill\ifx\theshortauthors\relax 
\theauthors\else\theshortauthors\fi\fi\fi}}\def\@evenhead{\@oddhead}
\def\@oddfoot{\small\lfoot\ifnum\count0=\startpage\copyright\ \gtp\hfill\else
\gt, Volume \thevolumenumber\ (\thevolumeyear)\hfill\fi}
\def\@evenfoot{\@oddfoot}
\makeatother
\fi

\newwrite\gtoutfile
\long\gdef\makeheadfile{  
{\def\\{, }\def\s{ }
\immediate\openout\gtoutfile head.xxx
\immediate\write\gtoutfile{Proxy-for: \ifx\theasciiauthors\relax
\theauthors\else\theasciiauthors\fi\s<\ifx\theasciiemail\relax\theemail\else\theasciiemail\fi>}
\immediate\write\gtoutfile{\noexpand\\}
\immediate\write\gtoutfile{Authors: \ifx\theasciiauthors\relax
\theauthors\else\theasciiauthors\fi}
{\def\\{ }\immediate\write\gtoutfile{Title: \ifx\theasciititle\relax
\thetitle\else\theasciititle\fi}}
\immediate\write\gtoutfile{Subj-class: GT or SG or MG etc}
\immediate\write\gtoutfile{MSC-class: \theprimaryclass\ifx\thesecondaryclass\relax\else, \thesecondaryclass\fi}
\immediate\write\gtoutfile{Journal-ref: Geom. Topol. \thevolumenumber
(\thevolumeyear) \startpage-\finishpage}
\immediate\write\gtoutfile{Comments: Published by Geometry and Topology at}
\immediate\write\gtoutfile{\s\s http://www.maths.warwick.ac.uk/gt/GTVol\thevolumenumber/paper\thepapernumber.abs.html}
\immediate\write\gtoutfile{\noexpand\\}
\immediate\write\gtoutfile{}
\ifx\theasciiabstract\relax
\immediate\write\gtoutfile{\theabstract}\else
\immediate\write\gtoutfile{\theasciiabstract}\fi
\immediate\write\gtoutfile{}
\immediate\write\gtoutfile{\noexpand\\}
\immediate\write\gtoutfile{}
\immediate\closeout\gtoutfile}}  

\def\maketitlepage{\maketitlep\makeheadfile}
\let\maketitle\maketitlepage